\newif\ifsubmit
\submitfalse

\ifsubmit
    \RequirePackage{fix-cm}
    \documentclass[smallextended, hidelinks]{svjour3}       
    \smartqed  
\else
    \documentclass[12pt, hidelinks]{article}
\fi

\usepackage{hyperref}
\usepackage{mathtools}
\usepackage{adjustbox}
\ifsubmit\else
    \usepackage[utf8]{inputenc}
    \usepackage[T1]{fontenc}
    \usepackage{fullpage}
\fi

\ifsubmit
\else
    \usepackage[style=alphabetic,backend=bibtex,maxbibnames=99,maxalphanames=99]{biblatex}
    \bibliography{refs}
\fi

\usepackage{amsmath, amssymb, url, float, float}

\usepackage[small,bf]{caption}
\usepackage{subcaption}
\usepackage{graphicx, xcolor}
\usepackage{tikz, pgfplots}
\pgfplotsset{
    compat=1.16,
}
\usetikzlibrary{topaths,intersections,calc}
\usetikzlibrary{positioning,chains,fit,shapes,calc,arrows.meta}

\let\phi\varphi
\let\eps\varepsilon

\newcommand{\BEAS}{\begin{eqnarray*}}
    \newcommand{\EEAS}{\end{eqnarray*}}
    \newcommand{\BEA}{\begin{eqnarray}}
    \newcommand{\EEA}{\end{eqnarray}}
    \newcommand{\BEQ}{\begin{equation}}
    \newcommand{\EEQ}{\end{equation}}
    \newcommand{\BIT}{\begin{itemize}}
    \newcommand{\EIT}{\end{itemize}}
    \newcommand{\BNUM}{\begin{enumerate}}
    \newcommand{\ENUM}{\end{enumerate}}
    
    \newcommand{\cf}{{\it cf.}}
    \newcommand{\eg}{{\it e.g.}}
    \newcommand{\ie}{{\it i.e.}}

    \newcommand{\ones}{\mathbf 1}
    
    \newcommand{\reals}{\mathbf{R}}

    \newcommand{\bmat}[1]{\begin{bmatrix}#1\end{bmatrix}}
    
    \newcommand{\dom}{\mathop{\bf dom}}

    \newcommand{\relint}{\mathop{\bf rel int}}

    \newif\iftodos

    \newcommand{\abs}[1]{\lvert{#1}\rvert}

\title{Convex Network Flows}
\ifsubmit
    \author{
        Theo Diamandis \and 
        Guillermo Angeris \and 
        Alan Edelman
    }
\else
    \author{
        Theo Diamandis \\ \texttt{\small tdiamand@mit.edu} \and 
        Guillermo Angeris \\ \texttt{\small gangeris@baincapital.com} \and 
        Alan Edelman \\ \texttt{\small edelman@mit.edu}
    }
\fi

\ifsubmit
    \institute{
        T. Diamandis \at MIT \email{tdiamand@mit.edu} \and
        G. Angeris \at Bain Capital \and
        A. Edelman \at MIT
    }
    \date{Received: date / Accepted: date}
\else
    \date{March 2024}
\fi

\begin{document} 
\maketitle 

\begin{abstract}
    We introduce a general framework for flow problems over hypergraphs. In our
    problem formulation, which we call the \emph{convex flow problem}, we have a
    concave utility function for the net flow at every node and a concave
    utility function for each edge flow. The objective is to maximize the sum of
    these utilities, subject to constraints on the flows allowed at each edge,
    which we only assume to be a convex set. This framework not only includes
    many classic problems in network optimization, such as max flow, min-cost
    flow, and multi-commodity flows, but also generalizes these problems to
    allow, for example, concave edge gain functions. In addition, our
    framework includes applications spanning a number of fields:
    optimal power flow over lossy networks, routing and resource allocation in
    ad-hoc wireless networks, Arrow-Debreu Nash bargaining, and order routing through
    financial exchanges, among others. We show that the convex flow problem
    has a dual with a number of interesting interpretations, and that this
    dual decomposes over the edges of the hypergraph. Using this decomposition, 
    we propose a fast solution algorithm that parallelizes over the edges and admits a
    clean problem interface. We provide an open source implementation of this
    algorithm in the Julia programming language, which we show is significantly
    faster than the state-of-the-art commercial convex solver Mosek.

    \ifsubmit
        \keywords{
            Network flows \and 
            Convex optimization \and 
            Dual decomposition \and 
            Generalized flows \and 
            Hypergraphs
        }
        \subclass{90C35 \and 90C25 \and 90C46 \and 05C21}
    \fi
\end{abstract}

\section{Introduction}
Network flow models describe a wide variety of common scenarios in computer
science, operations research, and other fields: from routing trucks to routing
bits. An extensive literature has developed theory, algorithms, and
applications for the case of linear flows over graphs. (See,
\eg,~\cite{ahuja1988network}, \cite{williamson2019network}, and references
therein.) However, the modeling capability of linear network flows is
significantly limited. For example, in many applications, the marginal flow out
of an edge decreases as the flow into this edge increases; \ie, the output from
the edge, as a function of its input, is concave. This property can be observed
in physical systems, such as power networks, where increasing the power through
a transmission line increases the line's loss, and in economic systems, such as
financial markets, where buying more of an asset increases the price of that
asset, resulting in a worse exchange rate. Additionally, there are many
applications where the flows through multiple edges connected to a single node
are nonlinearly related. For example, in a wireless network, a transmitter has
a power constraint across all of its links. Alternatively, in economics,
utilities may be superadditive when goods are complements. The fact that
classical network flow models cannot incorporate these well-studied
applications suggests that there is a natural generalization that can.

In this paper, we introduce the \emph{convex flow problem}, a generalization of
the network flow problem that significantly expands its modeling power. We
introduce two key ideas which, taken together, allow our framework to model many
additional problems present in the literature. First, instead of a graph, we
consider flows over a hypergraph---a graph where an edge may connect more than
two vertices. Second, we consider the allowable flows for each edge (which may
contain more than two vertices) to be a general convex set. This setup includes,
as special cases, the linear relationship studied in most network flow problems
and the concave, monotonic increasing edge input-output functions studied
in~\cite{shigeno2006maximum,vegh2014concave}. Our framework also encompasses a
number of other problems in networked physical and economic systems previously
studied in the literature. In many cases, it offers immediate generalizations or
more succinct formulations. We outline examples from a number of fields,
including power systems, wireless networks, Fisher markets, and financial asset
networks.

The convex flow problem we introduce is a convex optimization problem which
can, in practice, be efficiently solved. Our framework preserves the overall
network structure present in the problem and provides several interesting
insights. These insights, in turn, allow us to develop an efficient algorithm
for solving the convex flow problem. We show that the dual problem decomposes
over the network's edges, which leads to a highly-parallelizable algorithm that
can be decentralized. Importantly, this algorithm has a clean problem
interface: we only need access to (1) a Fenchel conjugate-like function of the
objective terms and (2) the solution to a simple subproblem for each edge.
These subproblem evaluations can be parallelized and have
efficiently-computable (and often closed form) solutions in many
applications. As a result, our algorithm enjoys better scaling and
order-of-magnitude faster solve times than commercial solvers like Mosek.

\paragraph{Outline.} We introduce a general framework for optimizing convex
flows over hypergraph structures, where each edge may connect more than two
vertices, in section~\ref{sec:the-problem}. In section~\ref{sec:apps}, we show
that this framework encompasses a number of problems previously studied in the
literature such as minimum cost flow and routing in wireless networks, and, in
some cases, offers immediate generalizations. We find a specific dual problem in
section~\ref{sec:dual-problem}, which we show has many useful interpretations
and decomposes nicely over the edges of the network. In
section~\ref{sec:algorithm}, we introduce an efficient algorithm that makes use
of this decomposition. This algorithm includes an efficient method to handle
edges which connect only two nodes, along with a method to recover a solution to
the original problem, using the solution to the dual, when the problem is not
strictly convex. Finally, we conclude with some numerical examples in
section~\ref{sec:numerical-examples}. This paper is accompanied by an open
source implementation of the solver in the Julia programming language.

\subsection{Related work}
The classic linear network flow problem has been studied extensively and we
refer the reader to~\cite{ahuja1988network} and~\cite{williamson2019network}
for a thorough treatment. In the classic setting, edges connecting more than
two vertices can be modeled by simply augmenting the graph with additional
nodes and two-node edges. While nonlinear cost functions have also been
extensively explored in the literature (\eg, see~\cite{bertsekas1998network}
and references therein), nonlinear edge flows---when the flow out of an edge is
a nonlinear function of the flow into it---has received considerably less
attention despite its increased modeling capability.

\paragraph{Nonlinear edge flows.} Extending the network flow problem to include
nonlinear edge flows was first considered by
Truemper~\cite{truemper1978optimal}. Still, work in the subsequent decades
mainly focused on the linear case---when the flow out of an edge is a linear
function of the flow into that edge---possibly with a convex cost function in
the objective. (See, for example, \cite{bertsekas1998network} and references
therein.) More recently, Shigeno~\cite{shigeno2006maximum} and
V\'{e}gh~\cite{vegh2014concave} considered the maximum flow problem where the flow
leaving and edge is a concave function of the flow entering that edge and
proposed theoretically efficient algorithms tailored to this case. This problem
is a special case of the convex flow problem we introduce in this work.
The nonlinear network flow problem has also appeared in a number of
applications, which we refer to in the relevant sections.

\paragraph{Dual decomposition methods for network flows.} The use of dual
decomposition methods for network flow problems has a long and rich history,
dating back to Kuhn's `Hungarian method' for the assignment
problem~\cite{kuhn1955hungarian}. The optimization community has explored these
methods extensively for network optimization problems (\eg, see~\cite[\S6,
\S9]{bertsekas1998network}) and, more generally, for convex optimization
problems with coupling constraints (\eg,
see~\cite[\S7]{bertsekasNonlinearProgramming2016}). These methods have also
been applied to many network problems in practice. For example, they have
facilitated the analysis and design of networking protocols, such as those used
for TCP congestion control~\cite{chiang2007layering}.
These protocols are, in essence, distributed, decentralized algorithms for
solving some global optimization problem. 

\paragraph{Extended monotropic programming.} Perhaps most related to our
framework is the extended monotropic programming problem, introduced by
Bertsekas~\cite{bertsekas2008extended}, of which our convex flow
problem is a special case. Both the convex flow problem and the
extended monotropic programming problem generalize Rockafellar's monotropic
programming problem~\cite{rockafellar1984network}. The strong duality result
of~\cite{bertsekas2008extended}, therefore, applies to our convex flow problem
as well, and we make this connection explicit in
appendix~\ref{app:extended-monotropic}. Although the convex flow problem we
introduce is a special case of the extended monotropic programming problem, our
work differs from that of Bertsekas along a number of dimensions. First, we
construct a different dual optimization problem which has a number of nice
properties. Second, this dual leads to a different algorithm than the one
developed in~\cite{bertsekas2008extended} and~\cite[\S4]{bertsekas2015convex},
and our dual admits an easier-to-implement interface with simpler
`subproblems'. Finally, while the application to multi-commodity flows is
mentioned in~\cite{bertsekas2008extended}, we show that our framework
encompasses a number of other problems in networked physical and economic
systems previously studied in the literature, and we numerically illustrate the
benefit of our approach.


\section{The convex flow problem}\label{sec:the-problem}
In this section, we introduce the convex flow problem, which 
generalizes a number of classic optimization problems in graph theory, including 
the maximum flow problem, the minimum cost flow problem, the multi-commodity 
flow problem, and the monotropic programming problem, among others. Our 
generalization builds on two key ideas: first, instead of a graph, we consider a 
hypergraph, where each edge can connect more than two nodes, and, second, we 
represent the set of allowable flows for each edge as a convex set. These two 
ideas together allow us to model many practical applications which have 
nonlinear relationships between flows.

\begin{figure}[h]
    \centering
    \hfill
    \adjustbox{max width=0.48\textwidth}{
        \begin{tikzpicture}
            \node (v1) at (0,2) {};
            \node (v2) at (3,2.5) {};
            \node (v3) at (0,0) {};
            \node (v4) at (4,-0.5) {};
        
            \begin{scope}[fill opacity=0.8]
            \filldraw[fill=blue!70] ($(v1)+(0.5,1)$) 
                to[out=0,in=180] ($(v3) + (1.5,0)$)
                to[out=0,in=90] ($(v4) + (1,-0.5)$)
                to[out=270,in=0] ($0.5*(v3) + 0.5*(v4) + (1,-1)$)
                to[out=180,in=270] ($(v3) + (-1,0)$)
                to[out=90,in=180] ($(v1)+(0.5,1)$);
            \node at (2,-0.5) {\footnotesize $e_1$};
            \filldraw[fill=red!70] ($(v1)+(0,0.5)$)
                to[out=0,in=90] ($.5*(v1)+.5*(v3) + (0.55,-0.25)$)
                to[out=270,in=0] ($(v3) + (0,-1)$)
                to[out=180,in=270] ($.5*(v1)+.5*(v3) + (-0.55,-0.25)$)
                to[out=90,in=180] ($(v1)+(0,0.5)$);
            \node at (0,1) {\footnotesize $e_2$};
            \filldraw[fill=green!70] ($(v2) + (-.25,1)$)
                to[out=0,in=90] ($.5*(v2)+.5*(v4) + (0.5,0)$)
                to[out=270,in=0] ($(v4) + (0.5,-.85)$)
                to[out=180,in=270] ($.5*(v2)+.5*(v4) + (-0.75,0)$)
                to[out=90,in=180] ($(v2) + (-.25,1)$);
            \node at ($.5*(v2)+.5*(v4)$) {\footnotesize $e_3$};
            \end{scope}
        
            \foreach \v in {1,2,3,4} {
                \fill (v\v) circle (0.1);
            }
        
            \fill (v1) circle (0.05) node [below] {$v_1$};
            \fill (v2) circle (0.05) node [above] {$v_2$};
            \fill (v3) circle (0.05) node [below] {$v_3$};
            \fill (v4) circle (0.05) node [below] {$v_4$};
        
        \end{tikzpicture}
    }
    \hfill
    \adjustbox{max width=0.48\textwidth}{
        
        \begin{tikzpicture}[thick,
            cfmm/.style={draw,circle},
            every fit/.style={ellipse,draw,inner sep=-2pt,text width=2cm},
            shorten >= 3pt,shorten <= 3pt,
          ]
          
          \node at (0,2.25) {edges};
          \begin{scope}[yshift=+10mm, start chain=going below, node distance=5mm]
            \foreach \i in {1,2,3}
                \node[cfmm,on chain] (e\i) {};
            \node at ($(e1) + (-0.5,0)$) [left] {\footnotesize $e_1$};
            \node at ($(e2) + (-0.5,0)$) [left] {\footnotesize $e_2$};
            \node at ($(e3) + (-0.5,0)$) [left] {\footnotesize $e_3$};
          \end{scope}
          
          \node at (4,2.25) {nodes};
          \begin{scope}[xshift=4cm,yshift=+12.5mm,start chain=going below,node distance=5mm]
            \foreach \i in {1,2,...,4}
                \node[on chain, circle, draw, fill=black] (v\i) {};
            \node at ($(v1) + (0.5,0)$) [right] {$v_1$};
            \node at ($(v2) + (0.5,0)$) [right] {$v_2$};
            \node at ($(v3) + (0.5,0)$) [right] {$v_3$};
            \node at ($(v4) + (0.5,0)$) [right] {$v_4$};
          \end{scope}
        
          
          
          \draw[blue] (e1) -- (v1);
          \draw[blue] (e1) -- (v3);
          \draw[blue] (e1) -- (v4);
          
          \draw[red] (e2) -- (v1);
          \draw[red] (e2) -- (v3);
          
          \draw[green] (e3) -- (v2);
          \draw[green] (e3) -- (v4);

          \node at (0,-2) {};
          
        \end{tikzpicture}
    }
    \hfill
    \caption{A hypergraph with 4 nodes and 3 edges (left) and its
    corresponding bipartite graph representation (right).}
    \label{fig:hypergraph}
\end{figure}
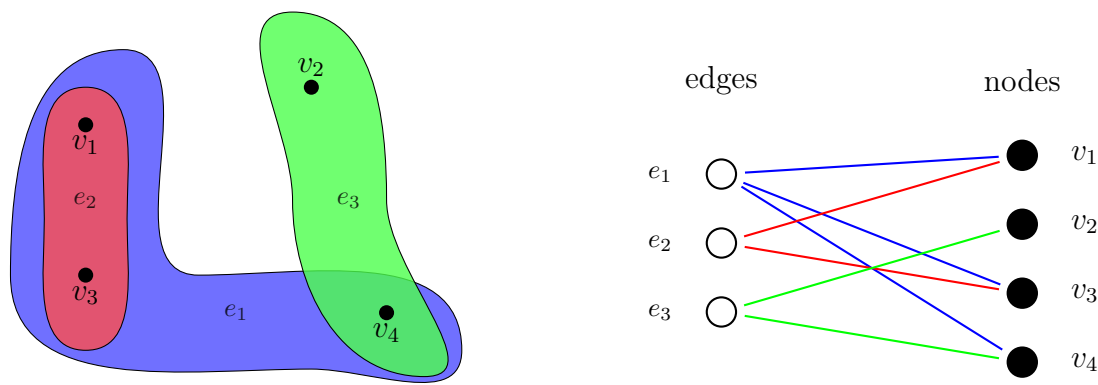

\paragraph{Hypergraphs.} We consider a hypergraph with $n$ nodes and $m$
hyperedges. Each hyperedge (which we will refer to simply as an `edge' from
here on out) connects some subset of the $n$ nodes.  This hypergraph may also
be represented as a bipartite graph with $n + m$ vertices, where the first
independent set contains $n$ vertices, each corresponding to one of the $n$
nodes in the hypergraph, and the second independent set contains the
remaining $m$ vertices, corresponding to the $m$ edges in the hypergraph. An
edge in the bipartite graph exists between vertex $i$, in the first independent
set, and vertex $j$, in the second independent set, if, and only if, in the
corresponding hypergraph, node $i$ is incident to (hyper)edge $j$.
Figure~\ref{fig:hypergraph} illustrates these two representations. From the
bipartite graph representation, we can easily see that the labeling of `nodes'
and `edges' in the hypergraph is arbitrary, and we will sometimes switch these
labels based on convention in the applications. While this section presents the
bipartite graph representation as a useful perspective for readers, it is
not used in what follows.

\paragraph{Flows.} On each of the edges in the graph, $i = 1, \dots, m$, we
denote the \emph{flow} across edge $i$ by a vector $x_i \in \reals^{n_i}$,
where $n_i \ge 2$ is the number of nodes incident to edge $i$. Each of these
edges $i$ also has an associated closed, convex set $T_i \subseteq
\reals^{n_i}$, which we call the \emph{allowable flows} over edge $i$, such
that only flows $x_i \in T_i$ are feasible. (We often also have that $0 \in T_i$,
\ie, we have the option to not use an edge.) By convention, we will use
positive numbers to denote flow out of an edge (equivalently, into a node) and
negative numbers to denote flow into an edge (equivalently, out of a node). For
example, in a standard graph, every edge connects exactly two vertices, so $n_i
= 2$. If $1$ unit of flow travels from the first node to the second node
through an edge, the flow vector across that edge is
\[
    x_i = \bmat{-1 \\ 1}.
\] 
If this edge is bidirectional (\ie, if flow can travel in either direction),
lossless, and has some upper bound $b_i > 0$ on the flow (sometimes called the
`capacity'), then its allowable flows are
\[
    T_i = \{z \in \reals^2 \mid z \le b_i\ones ~~ \text{and} ~~ z_1 + z_2 = 0\}.
\]
While this formalism may feel overly cumbersome when dealing with standard
graphs, it will be useful for working with hypergraphs.

\paragraph{Local and global indexing.}
We denote the number of nodes incident to edge $i$ by $n_i$. This set of `local'
incident nodes is a subset of the `global' set of $n$ nodes in the hypergraph. 
To connect the local node indices to the global node indices, we introduce 
matrices $A_i \in \reals^{n \times n_i}$. In particular, we define $(A_i)_{jk} = 1$ 
if node $j$ in the global index corresponds to node $k$ in the local index,
and $(A_i)_{j k} = 0$, otherwise.
For example, consider a hypergraph with 3 nodes. If edge $i$ connects nodes 2
and 3, then
\[
    A_i = \bmat{0 & 0 \\ 1 & 0 \\ 0 & 1 } = \bmat{ \vline & \vline \\e_2 & e_3 \\ \vline & \vline}.
\]
Written another way, if the $k$th node in the edge corresponds to global
node index $j$, then the $k$th column of $A_i$, is the $j$th unit basis vector,
$e_j$. Note that the ordering of nodes in the local indices need not be the 
same as the global ordering.

\paragraph{Net flows.}
By summing the flow in each edge, after mapping these flows to the
global indices, we obtain the \emph{net flow vector}
\[
    y = \sum_{i=1}^m A_ix_i.
\]
We can interpret $y$ as the netted flow across the hypergraph. If $y_j > 0$,
then node $j$ ends up with flow coming into it. (These nodes are often called 
\emph{sinks}.) Similarly, if $y_j < 0$, then node $j$ must provide some flow to 
the network. (These nodes are often called \emph{sources}.) Note that a node $j$
with $y_j = 0$ may still have flow passing through it; zero net flow only 
means that this node is neither a source nor a sink.

\paragraph{Utilities.} Now that we have defined the individual edge flows $x_i$
and the net flow vector $y$, we introduce utility functions for each.
First, we denote the \emph{network utility} by $U : \reals^n \to \reals
\cup \{-\infty\}$, which maps the net flow vector $y$ to a utility value,
$U(y)$. Infinite values denote constraints: any flow with $U(y) = -\infty$ is
unacceptable. We also introduce a utility function for each edge, $V_i :
\reals^{n_i} \to \reals \cup \{-\infty\}$, which maps the flow $x_i$ on edge
$i$ to a utility, $V_i(x_i)$. We require that both $U$ and the $V_i$ are 
concave, nondecreasing functions. This restriction is not as strong as it may
seem; we may also minimize convex nondecreasing cost 
functions with this framework.

\paragraph{Convex flow problem.} The \emph{convex flow problem}
seeks to maximize the sum of the network utility and the
individual edge utilities, subject to the constraints on the allowable flows:
\begin{equation}
    \label{eq:main-problem}
    \begin{aligned}
        & \text{maximize} && {\textstyle U(y) + \sum_{i=1}^m V_i(x_i)}\\
        & \text{subject to} && {\textstyle y = \sum_{i=1}^m A_i x_i}\\
        &&& x_i \in T_i, \quad i=1, \dots, m.
    \end{aligned}
\end{equation}
Here, the variables are the edge flows $x_i \in \reals^{n_i}$, for $i = 1,
\dots, m$, and the net flows $y \in \reals^n$. Each of these edges can be
thought of as a subsystem with its own local utility function $V_i$. The
individual edge flows $x_i$ are local variables, specific to the $i$th
subsystem. The overall system, on the other hand, has a utility that is a
function of the net flows $y$, the global variable. As we will see in what
follows, this structure naturally appears in many applications and lends itself
nicely to parallelizable algorithms. Note that, because the
objective is nondecreasing in all of its variables,
a solution $\{x^\star_i\}$ to problem~\eqref{eq:main-problem} will almost always have
$x_i^\star$ at the boundary of the feasible flow set $T_i$. If an $x_i^\star$
were in the interior, we could increase its entries without decreasing the
objective value until we hit the boundary of the corresponding $T_i$, assuming
some basic conditions on $T_i$ (\eg, $T_i$ does not contain a strictly
positive ray).


\section{Applications}\label{sec:apps}
In this section, we give a number of applications of the convex flow
problem~\eqref{eq:main-problem}. We first show that many classic
optimization problems in graph theory are special cases of this problem. Then, 
we show that the convex flow problem models problems in a variety of 
fields including power systems, communications, economics, and finance, among 
others. We start with simple special cases and gradually build up to those that 
are firmly outside the traditional network flows literature.

\subsection{Maximum flow and friends}\label{sec:app-max-flow}
In this subsection, we show that many classic network flow problems are special
cases of problem~\eqref{eq:main-problem}. We begin with a standard setup that
will be used for the rest of this subsection.

\paragraph{Edge flows.}
We consider a directed graph with $m$ edges and $n$ nodes, which we assume to
be connected. Recall that we denote the flow over edge $i$ by the vector $x_i
\in \reals^2$. We assume that edge $i$'s flow has upper bound $b_i \ge 0$, so
the set of allowable flows is
\begin{equation}\label{eq:Ti-simple-graph}    
    T_i = \{
        z \in \reals^2 \mid 0 \le z_2 \le b_i ~~ \text{and} ~~ z_1 + z_2 = 0
    \}.
\end{equation}
With this framework, it is easy to see how gain factors or other
transformations can be easily incorporated into the problem. For example, we
can instead require that $\alpha z_1 + z_2 = 0$ where $\alpha > 0$ is some gain
or loss factor. Note that if the graph is instead undirected, with each set of
two directed edges replaced by one undirected edge, the allowable flows for
each pair of directed edges can be combined into the set
\[ 
    T_i = \{
        z \in \reals^2 \mid z \le b_i \ones ~~ \text{and} ~~ z_1 + z_2 = 0
    \},
\]
which is the Minkowski sum of the two allowable flows in the directed case, one
for each direction. For what follows, we only consider directed graphs, but the
extension to the undirected case is straightforward.

\paragraph{Net flow.}
To connect these edge flows to the net flow we use the
matrices $A_i \in \{0,1\}^{n \times 2}$ for each edge $i = 1, \dots, m$ such
that, if edge $i$ connects node $j$ to node $k$ (assuming the direction of the
edge is from node $j$ to node $k$), then we have
\begin{equation}\label{eq:Ai-simple-graph}
    A_i = \bmat{ \vline & \vline \\e_{j} & e_{k} \\ \vline & \vline}.
\end{equation}
Using these matrices, we write the net flow through the network as the sum of
the edge flows:
\[
    y = \sum_{i=1}^m A_i x_i.
\]

\paragraph{Conservation laws.} One important consequence of the definition
of the allowable flows $T_i$ is that there is a corresponding \emph{local conservation law}:
for any allowable flow $x_i \in T_i$, we have that
\[
    \ones^Tx_i = (x_i)_1 + (x_i)_2 = 0,
\]
by definition of the set $T_i$. Since the $A_i$ matrices are simply selector
matrices, we therefore have that $\ones^TA_ix_i = 0$ whenever $x_i \in T_i$,
which means that we can turn the local conservation law above into a
\emph{global conservation law}:
\begin{equation}\label{eq:global-conservation}
    \ones^Ty = \sum_{i=1}^m \ones^TA_ix_i = 0,
\end{equation}
where $y$ is the corresponding net flow, for any set of allowable flows $x_i
\in T_i$, for $i=1, \dots, m$. We will use this fact to show that feasible
flows are, indeed, flows through the network in the `usual' sense. Conversely,
we can find conservation laws for a given convex flow problem, which we discuss
in appendix~\ref{app:conservation}.

\subsubsection{Maximum flow}
Given a directed graph, the maximum flow problem seeks to find the maximum amount
of flow that can be sent from a designated source node to a designated sink node. The problem
can model many different situations, including transportation network routing, 
matching, and resource allocation, among others. It dates back to the work of 
Harris and Ross~\cite{harris1955fundamentals} to model Soviet railway networks 
in a report written for the US Air Force and declassified in 1999, at the 
request of Schrijver~\cite{schrijver2002history}. While well-known to be a 
linear program at the time~\cite{ford1956maximal} (and therefore solvable with 
the simplex method), specialized methods were quickly developed~\cite{ford1957simple}. 
The maximum flow problem has been extensively studied by the operations research 
and computer science communities since then.

\paragraph{Flow conservation.} Relabeling the graph such that the source node
is node $1$ and the sink node is node $n$, we write the net flow conservation
constraints as the set
\begin{equation}
    \label{eq:flow-conservation-constraints}
    S = \{y \in \reals^n \mid y_1 + y_n \ge 0, \quad y_j \ge 0 \quad \text{for all} ~~ j \ne 1, n\}.
\end{equation}
Note that this set $S$ is convex as it is the intersection of halfspaces (each
of which is convex), and its corresponding indicator function, written
\[
    I_S(y) = \begin{cases}
        0 & y \in S\\
        +\infty & \text{otherwise},
    \end{cases}
\]
is therefore also convex. This indicator function is nonincreasing in that, if
$y' \ge y$ then $I_S(y') \le I_S(y)$ by definition of the set $S$. Thus,
its negation, $-I_S$, is nondecreasing and concave.

\paragraph{Problem formulation.}
The network utility function in the maximum flow problem is to maximize the flow 
into the terminal node while respecting the flow conversation constraints:
\[
    U(y) = y_n - I_{S}(y).
\]
From the previous discussion, this utility function is concave and
nondecreasing. We set the edge utility functions to be zero, $V_i = 0$ for all
$i = 1, \dots, m$, to recover the maximum flow problem (see, for example, 
\cite[Example 1.3]{bertsekas1998network}) in our framework:
\begin{equation}\label{eq:max-flow}
    \begin{aligned}
        & \text{maximize} && y_n - I_S(y)\\
        & \text{subject to} && {\textstyle y = \sum_{i=1}^m A_i x_i}\\
        &&& x_i \in T_i, \quad i=1, \dots, m,
    \end{aligned}
\end{equation}
where the sets $\{T_i\}$ are the feasible flow sets~\eqref{eq:Ti-simple-graph} 
and may be either directed or undirected.

\paragraph{Problem properties.} Any feasible point (\ie, one that satisfies the
constraints and has finite objective value) is a flow in that the net flow
through node $j$ is zero ($y_j = 0$) for any node that is not the source,
$j=1$, or the sink, $j=n$. To see this, note that, for any $j\ne 1,n$, if $y
\in S$, then
\[
    \ones^Ty \ge y_1 + y_n + y_j \ge y_j \ge 0.
\]
The first and third inequalities follow from the fact that $y \in S$ means that
$y_j \ge 0$ for all $j \ne 1, n$; the second follows from the fact that $y_1 +
y_n \ge 0$ from the definition of $S$ as well. From the conservation
law~\eqref{eq:global-conservation}, we know that $\ones^Ty = 0$, so $y_j = 0$
for every node $j$ that is not the source nor the sink. Therefore, $y$ is a
flow in the `usual' sense. A similar proof shows that $-y_1 = y_n$; \ie,
the amount provided by the source is the amount dissipated by the sink. Since
we are maximizing the total amount dissipated $y_n$, subject to the provided
capacity and flow constraints, the problem above corresponds exactly to
the standard maximum flow problem.

\subsubsection{Minimum cost flow}
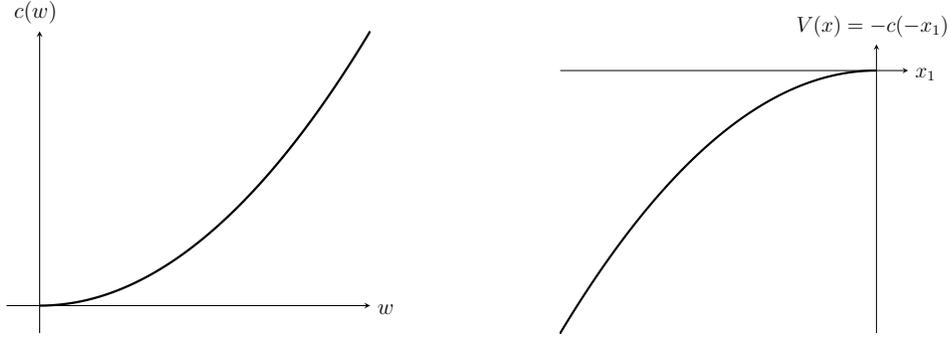
\begin{figure}[h]
    \centering
    \hfill
    \adjustbox{max width=0.33\textwidth}{
        \begin{tikzpicture}[scale=1.0]
    \begin{axis}[
        axis on top=true,
        xmin=-0.1, xmax=1,
        ymin=-0.1, ymax=1,
        axis lines=center,
        xlabel={$w$},
        xlabel style={at={(ticklabel cs:1)}, anchor=west},
        ylabel={$c(w)$},
        ylabel style={at={(ticklabel cs:1)}, anchor=south},
        grid=major,
        legend pos=outer north east,
        samples=200,
        domain=-4:4,
        xtick=\empty,
        ytick=\empty
        ]
        \addplot[black, very thick, domain=0:1] {(x)^2};
    \end{axis}
\end{tikzpicture}
    }
    \hfill
    \adjustbox{max width=0.33\textwidth}{
        \begin{tikzpicture}[scale=1.0]
    \begin{axis}[
        axis on top=true,
        xmin=-1, xmax=0.1,
        ymin=-1, ymax=0.1,
        axis lines=center,
        xlabel={$x_1$},
        xlabel style={at={(ticklabel cs:1)}, anchor=west},
        ylabel={$V(x) = -c(-x_1)$},
        ylabel style={at={(ticklabel cs:1)}, anchor=south},
        grid=major,
        legend pos=outer north east,
        samples=200,
        domain=-4:4,
        xtick=\empty,
        ytick=\empty
        ]
        \addplot[black, very thick, domain=-4:0] {-(x)^2};
    \end{axis}
\end{tikzpicture}
    }
    \hfill \null
    \caption{
        An example convex nondecreasing cost function $c(w) = w^2$ for $w \ge 0$
        (left) and its corresponding concave, nondecreasing edge utility
        function $V$ (right).
    }
    \label{fig:min-cost-function}
\end{figure}

The minimum cost flow problem seeks to find the cheapest way to route a given
amount of flow between specified source and sink nodes. We consider the same
setup as above, but with two modifications: first, we fix the value of the flow
from node $1$ to node $n$ to be at least some value $v \ge 0$; and second, we
introduce a convex, nondecreasing cost function for each edge $i$, denoted $c_i
:\reals_+ \to \reals_+$, which maps the flow on this edge to a cost. We modify
the flow conservation constraints to be
(\cf,~\eqref{eq:flow-conservation-constraints})
\[
    \tilde S = \{y \mid y_n \ge v,\quad y_1 + y_n \ge 0,\quad y_j \ge 0 \quad \text{for all} ~~ j\ne 1, n\}.
\]
Much like the previous, the negative indicator of this set, $-I_{\tilde S}$, is 
a concave, nondecreasing function. We take the edge flow utility function $V_i$ 
to be
\[
    V_i(x_i) = -c_i(-(x_i)_1),
\]
which is a concave nondecreasing function of $x_i$. (Recall that $(x_i)_1 \le
0$. We provide an example in figure~\ref{fig:min-cost-function}.) Modifying the
network utility function to be the indicator over this new set $\tilde S$,
\[
    U(y) = -I_{\tilde S}(y),
\]
we recover the minimum cost flow problem in our framework:
\[
    \begin{aligned}
        & \text{maximize} && {\textstyle {-I_{\tilde S}(y)} + \sum_{i=1}^{m} -c_i(-(x_i)_1)} \\
        & \text{subject to} && {\textstyle y = \sum_{i=1}^m A_i x_i}\\
        &&& x_i \in T_i, \quad i=1, \dots, m.
    \end{aligned}
\]
Here, as before, the sets $\{T_i\}$ are the directed feasible flow sets defined
in~\eqref{eq:Ti-simple-graph}.

\subsubsection{Concave edge gains}
\begin{figure}[h]
    \centering
    \hfill
    \adjustbox{max width=0.33\textwidth}{
        \begin{tikzpicture}[scale=1.0]
    \begin{axis}[
        axis on top=true,
        xmin=-0.1, xmax=1,
        ymin=-0.1, ymax=1,
        axis lines=center,
        xlabel={$w$},
        ylabel={$\gamma(w)$},
        grid=major,
        legend pos=outer north east,
        samples=200,
        domain=-4:4,
        xtick=\empty,
        ytick=\empty
        ]
        \addplot[black, very thick, domain=0:1] {sqrt(x)};
    \end{axis}
\end{tikzpicture}
    }
    \hfill
    \adjustbox{max width=0.33\textwidth}{
        \begin{tikzpicture}[scale=1.0]
    \begin{axis}[
        axis on top=true,
        xmin=-1, xmax=0.1,
        ymin=-0.1, ymax=1,
        axis lines=center,
        xlabel={$z_1$},
        ylabel={$z_2$},
        grid=major,
        legend pos=outer north east,
        samples=200,
        domain=-4:4,
        xtick=\empty,
        ytick=\empty
        ]
        \addplot[black, very thick, domain=-1:0] {sqrt(-x)};
        \addplot[gray!30, domain=-1:0, fill, opacity=0.6, draw=none] {sqrt(-x)} -| (-1, 0) -- cycle;
        \node at (-0.6, 0.4) {$T$};

        \draw[black, very thick] (0, 0) -- (-1, 0);
    \end{axis}
\end{tikzpicture}
    }
    \hfill \null
    \caption{
        An example concave edge gain function $\gamma(w) = \sqrt{w}$ (left) and
        the corresponding allowable flows (right).
    }
    \label{fig:gain-function}
\end{figure}
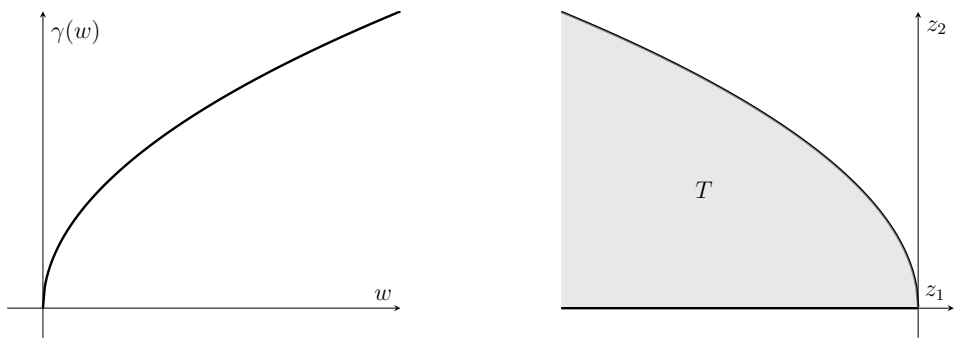

We can generalize the maximum flow problem and the minimum cost flow problem
to include concave, nondecreasing edge input-output functions, as
in~\cite{shigeno2006maximum,vegh2014concave}, by modifying the sets of feasible
flows. We denote the edge input-output functions by $\gamma_i
: \reals_+ \to \reals_+$. (For convenience, negative arguments to $\gamma_i$ are
equal to negative infinity.) If $w$ units of flow enter edge $i$, then
$\gamma_i(w)$ units of flow leave edge $i$. In this case, we can write the set
of allowable flows for each edge to be
\[
    T_i = \{ z \in \reals^2 \mid z_2 \le \gamma_i(-z_1) \}.
\]
We provide an example in figure~\ref{fig:gain-function}.
The inequality has the following interpretation: the magnitude of the flow 
out of edge $i$, given by $z_2 \ge 0$, can be any value not exceeding 
$\gamma_i(-z_1)$; however, we can `destroy' flow. From the problem properties 
presented in section~\ref{sec:the-problem}, there exists a solution such that 
this inequality is tight, since the utility function $U$ is nondecreasing. In 
other words, we can find a set of feasible flows $\{x_i\}$ such that
\[
    (x_i)_2 = \gamma_i(-(x_i)_1),
\]
for all edges $i = 1, \dots, m$.

\subsubsection{Multi-commodity flows}\label{sec:app-multi-commodity}
Thus far, all the flows have been of the same type. Here, we show that the
\emph{multi-commodity flow problem}, which seeks to route $K$ different
commodities through a network in an optimal way, is also a special case of the
convex flow problem. We denote the flow of these commodities over an edge $i$
by $x_i \in \reals^{2K}$, where the first $2$ elements denote the flow of the
first commodity through edge $i$, the next $2$ elements denote the flow of the
second, and so on. The set $T_i$ then allows us to specify joint constraints on
these flows. For example, we can model a total flow capacity by the set
\[
    T_i = \left\{x \in \reals^K ~\middle\vert~ \sum_{k=1}^K w_k x_{2k} \le b_i, ~ 0 \le x_{2k}, ~ \text{and} ~ x_{2k} = -x_{2k-1}~\text{for}~k=1, 2, \dots, K\right\},
\]
where $b_i$ denotes the capacity of edge $i$ and $w_k$ denotes the capacity
required per unit of commodity $k$. In other words, each commodity $k$ has a 
per-unit `weight' of $w_k$, and the maximum weighted capacity through edge $i$ 
is $b_i$. If $K=1$, then this set of allowable flows reduces to the original definition~\eqref{eq:Ti-simple-graph}. 
Additionally, note that $T_i$ is still a polyhedral set, but more complicated 
convex constraints may be added as well.

We denote the net flows at each node by a vector $y \in \reals^{nK}$, where the
first $K$ elements denote the net flows of the first commodity, the next $K$
elements denote the net flows of the second, and so on, while the $A_i$
matrices map the local flows of each good to the corresponding indices in $y$.
For example, if edge $i$ connects node $j$ to node $k$, the edge would have the
associated matrix $A_i \in \reals^{nK \times 2K}$ given by
\[
    A_i = \bmat{
        \vline & \vline & \vline & \vline && \vline & \vline \\
        e_{j} & e_{k} & e_{j+n} & e_{k+n} &\cdots& e_{j + (K-1)n} & e_{k + (K-1)n} \\
        \vline & \vline & \vline & \vline && \vline & \vline
    }.
\]
The problem is now analogous to those in the previous sections, only with $y$ 
and $x_i$ having larger dimension and $T_i$ modified as described above.

\subsection{Optimal power flow}\label{sec:app-opf}
The optimal power flow problem~\cite{wood2013power} seeks a cost-minimizing plan 
to generate power satisfying demand in each region. We consider a network of $m$ 
transmission lines (edges) between $n$ regions (nodes). We assume that the
region-transmission line graph is directed for simplicity.

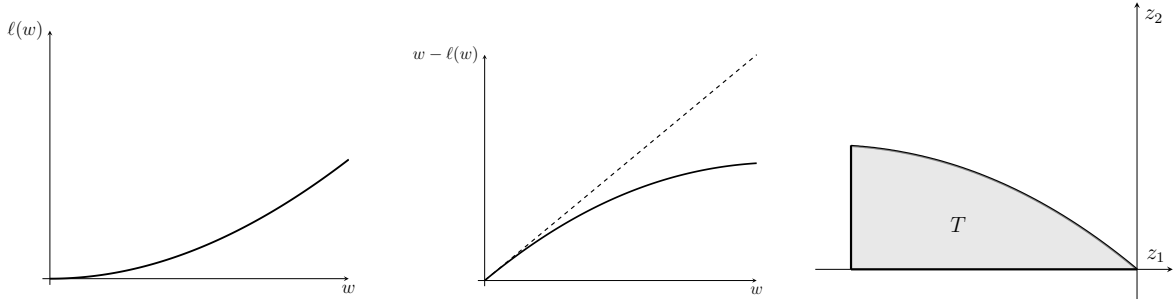
\begin{figure}
    \centering
    \hfill
    \adjustbox{max width=0.30\textwidth}{
        \input figures/opf-loss.tex
    }
    \hfill
    \adjustbox{max width=0.30\textwidth}{
        \input figures/opf-gain.tex
    }
    \hfill
    \adjustbox{max width=0.30\textwidth}{
        \begin{tikzpicture}[scale=1.0]
    \begin{axis}[
        axis on top=true,
        xmin=-4.5, xmax=0.5,
        ymin=-0.5, ymax=4.5,
        axis lines=center,
        xlabel={$z_1$},
        ylabel={$z_2$},
        grid=major,
        legend pos=outer north east,
        samples=200,
        domain=-4:4,
        xtick=\empty,
        ytick=\empty
        ]
        \addplot[black, very thick, domain=-4:0] {3*-x - 16*(ln(1 + exp(-x/4)) - ln(2))};
        \addplot[gray!30, domain=-4:0, fill, opacity=0.6, draw=none] {3*-x - 16*(ln(1 + exp(-x/4)) - ln(2))} -| (-4, 0) -- cycle;
        \node at (-2.5, 0.75) {$T$};

        \draw[black, very thick] (-4, 2.09) -- (-4, 0);
        \draw[black, very thick] (0, 0) -- (-4, 0);
    \end{axis}
\end{tikzpicture}
    }
    \hfill \null
    \caption{The power loss function (left), the corresponding power output (middle),
    and the corresponding set of allowable flows (right).}
    \label{fig:opf-loss}
\end{figure}

\paragraph{Line flows.} When power is transmitted along a line, the line heats
up and, as a result, dissipates power. As greater amounts of power are 
transmitted along this line, the line further heats up, which, in turn, causes 
it to dissipate even more power. We model this dissipation as a convex function 
of the power transmitted, which captures the fact that the dissipation increases 
as the power transmitted increases. We use the logarithmic power loss function 
from~\cite[\S 2.1.3]{stursberg2019mathematics}. With this loss function, the
`transport model' optimal power flow solution matches that of the more-complicated 
`DC model', assuming a uniform line material. (See~\cite[\S2]{stursberg2019mathematics}
for details and discussion.) The logarithmic loss function is given by
\[
    \ell_i(w) = \alpha_i \left(\log(1 + \exp(\beta_i w)) - \log 2\right) - 2w,
\]
where $\alpha_i$ and $\beta_i$ are known constants and $\alpha_i \beta_i = 4$ 
for each line $i$. This function can be easily verified to be convex,
increasing, and have $\ell_i(0) = 0$. The power output of a line with input $w$
can then be written as $w - \ell(w)$. We also introduce capacity constraints
for each line $i$, given by $b_i$. Taken together, for a given line $i$, the
power flow $x_i$ must lie within the set
\begin{equation}\label{eq:app-opf-Ti}
    T_i = \left\{ 
    z \in \reals^2 \mid 
    -b_i \le z_1 \le 0
    ~\text{and}~
    z_2 \le -z_1 - \ell_i(-z_1)
    \right\}, 
    \qquad i = 1, \dots, m.
\end{equation}
This set is convex, as it is the intersection of two halfspaces and the
epigraph of a convex function. Note that we relaxed the
line flow constraint to an inequality. This inequality has the following
physical interpretation: we may dissipate additional power along the line (for
example, by adding a resistive load), but in general we expect this inequality to
hold with equality, as discussed in~\S\ref{sec:the-problem}. Figure~\ref{fig:opf-loss}
shows a loss function and the corresponding edge's allowable flows.

\paragraph{Net flows.}
Each region $i = 1, \dots, n$ demands $d_i$ units of power.
In addition, region $i$ can generate power $p_i$ at cost $c_i : \reals \to \reals_+ \cup \{\infty\}$,
where infinite values denote constraints (\eg, a region may have a maximum power 
generation capacity). We assume that $c_i$ is convex and nondecreasing, with
$c_i(p_i) = 0$ for $p_i \le 0$ (\ie, we can dissipate power at zero cost). 
Similarly to the max flow problem in~\S\ref{sec:app-max-flow}, we take the
indexing matrices $A_i$ as defined in~\eqref{eq:Ai-simple-graph}.
To meet demand, we must have that
\[
    d = p + y, \qquad \text{where} \qquad  y =  \sum_{i=1}^m A_i x_i.
\]
In other words, the power produced, plus the net flow of power, must satisfy the 
demand in each region.
We write the network utility function as
\begin{equation}\label{eq:app-opf-U}
    U(y) = \sum_{i=1}^n -c_i(d_i - y_i).
\end{equation}
Since each $c_i$ is convex and nondecreasing, the utility function $U$ is concave and 
nondecreasing in $y$.
This problem can then be cast as a special case of the convex flow problem~\eqref{eq:main-problem}:
\[
\begin{aligned}
    & \text{maximize} && {\textstyle \sum_{i=1}^n -c_i(d_i - y_i)} \\
    & \text{subject to} && {\textstyle y = \sum_{i=1}^m A_i x_i}\\
    &&& x_i \in T_i, \quad i=1, \dots, m,
\end{aligned}
\]
with the same variables $\{x_i\}$ and $y$, zero edge utilities ($V_i =
0$), and the feasible flow sets $T_i$ given in~\eqref{eq:app-opf-Ti}.

\paragraph{Extensions and related problems.}
This model can be extended in a variety of ways.
For example, a region may have some joint capacity over all the power it outputs.
When there are constraints such as these resulting from the interactions between 
edge flows, a hypergraph model is more appropriate. On the other hand, simple 
two-edge concave flows model behavior in number of other types of networks: 
in queuing networks, throughput is a concave function of the input due to convex
delays~\cite[\S5.4]{bertsekas1992data}; similarly, in routing games~\cite[\S18]{roughgarden2007routing}, 
a convex cost function often implies a concave throughput; in perishable 
product supply chains, such as those for produce, increased volume leads to 
increased spoilage~\cite[\S2.3]{nagurney2022spatial}; and in reservoir networks~\cite[\S8.1]{bertsekas1998network},
seepage may increase as volume increases. Our framework not only can model these 
problems, but also allows us to easily extend them to more complicated settings.

\subsection{Routing and recourse allocation in wireless networks}\label{sec:app-wireless}
In many applications, standard graph edges do not accurately capture
interactions between multiple flows coming from a single node---there
may be joint, possibly nonlinear, constraints on all the flows involving this
node. To represent these constraints in our problem, we make use of the fact
that an edge may connect more than two nodes in~\eqref{eq:main-problem}. In this 
section, we illustrate this structure through the problem of jointly optimizing 
the data flows and the power allocations for a wireless network, heavily 
inspired by the formulation of this problem in~\cite{xiao2004simultaneous}.

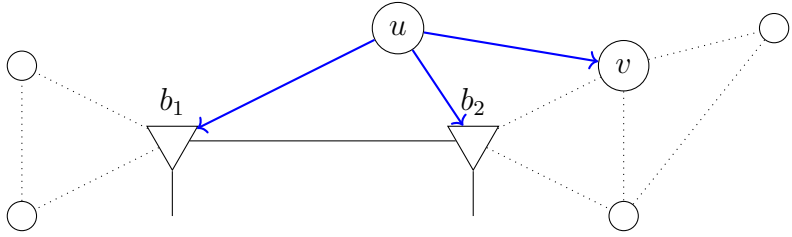
\begin{figure}
    \centering
    \adjustbox{max width=.7\textwidth}{
        \begin{tikzpicture}
    \node[draw, regular polygon, regular polygon sides=3, shape border rotate=180] (B1) at (0,0) {};
    \node[above] at (B1.north) {$b_1$};
    \draw (B1) -- ++(0, -1);
    \node[draw, regular polygon, regular polygon sides=3, shape border rotate=180] (B2) at (4,0) {};
    \node[above] at (B2.north) {$b_2$};
    \draw (B2) -- ++(0, -1);

    \node[draw, circle] (U1) at (-2,1) {};
    \node[draw, circle] (U2) at (-2,-1) {};
    \node[draw, circle] (U3) at (8,1.5) {};
    \node[draw, circle] (U4) at (6,1) {$v$};
    \node[draw, circle] (U5) at (6,-1) {};
    \node[draw, circle] (U6) at (3,1.5) {$u$};

    \draw[dotted] (B1) -- (U1);
    \draw[dotted] (U2) -- (U1);
    \draw[dotted] (B1) -- (U2);
    \draw[dotted] (U4) -- (U3);
    \draw[dotted] (U5) -- (U3);
    \draw[dotted] (B2) -- (U4);
    \draw[dotted] (B2) -- (U5);
    \draw[dotted] (U4) -- (U5);
    \draw[solid, <-, blue, thick] (U4) -- (U6);
    \draw[solid, <-, blue, thick] (B1) -- (U6);
    \draw[solid, <-, blue, thick] (B2) -- (U6);

    \draw (B1) -- (B2);
\end{tikzpicture}
    }
    \caption{Wireless ad-hoc network. The (outgoing) hyperedge associated with
    user $u$ is shown in blue, and the corresponding set of outgoing neighbors
    $O_u$ contains user $v$ and the two base stations, $b_1$ and $b_2$.}
    \label{fig:app-wireless-network}
\end{figure}

\paragraph{Data flows.}
We represent the topology of a data network by a directed graph with $n$ nodes
and $m = n$ edges: one for each node.
We want to route traffic from a particular source to a particular destination in
the network. (This model can be easily extended to handle multiple source-destination
pairs, potentially with different latency or bandwidth constraints, using the 
multi-commodity flow ideas discussed in~\S\ref{sec:app-multi-commodity}.)
We model the network with a hypergraph, where edge $i = 1, \dots, n$ is associated with 
node $i$ and connects $i$ to all its outgoing neighbors, which we denote by the
set $O_i$, as shown in figure~\ref{fig:app-wireless-network}. (In other words, if $j \in O_i$,
then node $j$ is a neighbor of node $i$.)
On each edge, we have a rate at which we transmit data, denoted by a vector 
$x_i \in \reals^{\abs{O_i} + 1}$, where the $k$th element of $x_i$ denotes the rate 
from node $i$ to its $k$th outgoing neighbor and the last component is the total
outgoing rate from node $i$. The net flow vector $y \in \reals^{n}$ can be written 
as
\[
    y = \sum_{i=1}^n A_i x_i,
\]
for indexing matrices $A_i$ in $\reals^{n \times (|O_i| + 1)}$, given by
\[
    A_i = \bmat{
        \vline &  & \vline & \vline \\
        e_{O_{i1}} & \dots & e_{O_{i|O_i|}} & e_{i} \\
        \vline &  & \vline & \vline
    },
\]
where $O_{ik}$ denotes the $k$th neighbor of $O_i$ in any order fixed
ahead of time.

\paragraph{Communications constraints.}
Hyperedges allow us to more easily model the communication constraints in the 
examples of~\cite{xiao2004simultaneous}. We associate some communication 
variables with each edge $i$. These variables might be, for example, power 
allocations, bandwidth allocations, or time slot allocations. We assume that the 
joint constraints on the transmission rate and the communication variables are 
some convex set. For example, take the communication variables to be 
$(p_i, w_i)$, where $p \in \reals^{\abs{O_i}}$ is a vector of power allocations 
and $w \in \reals^{\abs{O_i}}$ is a vector of bandwidth allocations to each 
of node $i$'s outgoing neighbors. We may have maximum power and bandwidth 
constraints, given by $p^{\mathrm{max}}_i$ and $w^{\mathrm{max}}_i$, so the set 
of feasible powers and bandwidths is
\[
    P_i = \{(p, w) \in \reals^{\abs{O_i}} \times \reals^{\abs{O_i}} \mid 
    ~p \ge 0,
    ~w \ge 0,
    ~\ones^Tp \le p^{\mathrm{max}},
    ~\ones^Tw \le w^{\mathrm{max}}
    \}.
\] 
These communication variables determine the rate at which node $i$ can transmit
data to its neighbors. For example, in the Gaussian broadcast channel with 
frequency division multiple access, this rate is governed by the Shannon 
capacity of a Gaussian channel~\cite{shannon1948mathematical}. The set of 
allowable flows can be written as
\[
    T_i = \left\{ (z, t) \in \reals^{\abs{O_i}} \times \reals ~\middle\vert~
        \ones^T z = -t, 
        ~~z \le w\circ \log_2\left(\ones + \frac{p}{\sigma w} \right),
        ~~
        (p, w) \in P_i
    \right\},
\]
where $\sigma \in \reals_+^n$ is a parameter that denotes the average power of
the noise in each channel, the logarithm and division, along with the
inequality, are applied elementwise, and $\circ$ denotes the elementwise
(Hadamard) product. The set $T_i$ is a convex set, as the logarithm is a concave
function and $w_k \log(1 + p_k/\sigma w_k)$, viewed as a function over the
$k$th element of each of the communication variables $(p, w)$, is the
perspective transformation of $\log(1 + p_k/\sigma)$, viewed as a function over
$p_k$, which preserves concavity~\cite[\S3.2.6]{cvxbook}. The remaining sets
are all affine or polyhedral, and intersections of convex sets are convex,
which gives the final result.

Importantly, the communication variables (here, the power allocations $p$ and
bandwidth allocations $w$) can be private to a node $i$; the optimizer only 
cares about the resulting public data flow rates $x_i \in T_i$. This structure 
not only simplifies the problem formulation but also hints at efficient, 
decentralized algorithms to solve this problem. We note that the hypergraph 
model allows us to also consider the general multicast case as well.

\paragraph{The optimization problem.}
Without loss of generality, denote the source node by $1$ and the sink node by
$n$. We may simply want to maximize the rate of data from the source to the
sink, in which case we can take the network utility function to be
\[
    U(y) = y_n - I_S(y),
\]
where the flow conversation constraints $S$ are the same as those of the
classic maximum flow problem, defined
in~\eqref{eq:flow-conservation-constraints}.
We may also use the functions $V_i$ to include utilities or costs associated with 
the transmission of data by node $i$. 
We can include communication variables in the objective as well by simply 
redefining the allowable flows $T_i$ to include the relevant communication 
variables and modifying the $A_i$'s accordingly to ignore these entries of $x_i$.
This modification is useful when we have costs associated with these 
variables---for example, costs on power consumption.
Equipped with the set of allowable flows and these utility functions, we can
write this problem as a convex flow problem~\eqref{eq:main-problem}.

\paragraph{Related problems.}
Many different choices of the objective function and constraint sets for
communication network allocation problems appear in the
literature~\cite{xiao2004simultaneous,bertsekas1998network}. This setup also
encompasses a number of other `resource allocation' problems where the network
structure isn't immediately obvious, one of which we discuss in the next
section.

\subsection{Market equilibrium and Nash bargaining}\label{sec:app-market}
Our framework includes and generalizes the concave network flow model used by 
V\'{e}gh~\cite{vegh2014concave} to study market equilibrium problems such as
Arrow-Debreu Nash bargaining~\cite{vazirani2012notion}. Consider a 
market with a set of $n_b$ buyers and $n_g$ goods. There is one divisible unit 
of each good to be sold. Buyer $i$ has a budget $b_i \ge 0$ and receives utility 
$u_{i}: \reals^{n_g}_+ \to \reals_+$ from some allocation $x_{i} \in [0,1]^{n_g}$ 
of goods. We assume that $u_{i}$ is concave and nondecreasing, with $u_i(0) = 0$ 
for each $i = 1, \dots, n_b$. An equilibrium solution to this market is an 
allocation of goods $x_{i} \in \reals^{n_g}$ for each buyer $i = 1, \dots, n_b$, 
and a price $p_j \in \reals_+$ for each good $j = 1, \dots, n_g$, such that: 
(1) all goods are sold; (2) all money of all buyers is spent; and (3) each buyer 
buys a `best' (\ie, utility-maximizing) bundle of goods, given these prices.

An equilibrium allocation for this market is given by a solution to the 
following convex program:
\begin{equation}
    \label{eq:app-market-equilibrium}
    \begin{aligned}
        &\text{maximize} && \sum_{i = 1}^{n_b} b_i \log(u_i(x_i)) \\
        &\text{subject to} && \sum_{i = 1}^{n_b} (x_i)_j = 1, \quad j = 1, \dots, n_g \\
        &&& x_i \ge 0, \quad i = 1, \dots, n_b.
    \end{aligned}
\end{equation}
Eisenberg and Gale~\cite{eisenberg1959consensus} proved that the optimality 
conditions of this convex optimization problem give the equilibrium conditions 
in the special case that $u_{i}(x)$ is linear (and therefore separable across 
goods), \ie,
\[
    u_i(x_i) = v_i^T x_i
\]
for constant weights $v_i \in \reals_+^{n_g}$. (We show the same result for the 
general case~\eqref{eq:app-market-equilibrium} in appendix~\ref{app:fisher-market}.)
The linear case is called the `linear Fisher market model'~\cite{vazirani2007combinatorial} 
and can be easily recognized as a special case of the standard maximum flow 
problem~\eqref{eq:max-flow} with nonnegative edge gain factors~\cite[\S3]{vegh2014concave}.

V\'{e}gh showed that all known extensions of the linear Fisher market model are
a special case of the generalized market problem~\eqref{eq:app-market-equilibrium}, 
where the utility functions $u_i$ are separable across the goods, \ie,
\[
    u_i(x_i) = \sum_{j = 1}^{n_g} u_{ij}((x_i)_j) 
\] 
for $u_{ij}: \reals \to \reals$ concave and increasing.
V\'{e}gh casts this problem as a maximum network flow problems with concave edge 
input-output functions.
We make a further extension here by allowing the utilities to be
concave nondecreasing functions of the entire basket of goods rather than sums
of functions of the individual allocations. This generalization allows us to 
model complementary goods and is an immediate consequence of our framework.

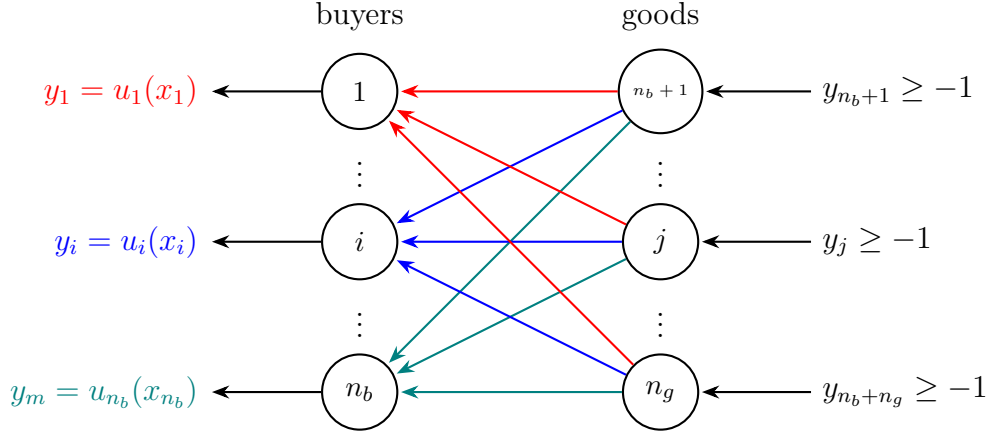
\begin{figure}
    \centering
        
        

        
        
        
        
        

\begin{tikzpicture}[
    ->,
    >=Stealth,
    shorten >=1pt,
    auto,
    node distance=2cm,
     thick,
    main node/.style={circle,
    draw,
    font=\sffamily\Large\bfseries}
]
    
    \node (buyers) at (0,0) {buyers};
    \node[main node, minimum size=1cm] (b1) at (0,-1) {\small $1$ };
    \node at (0,-2) {\vdots};
    \node[main node, minimum size=1cm] (bi) at (0,-3) {\small $i$ };
    \node at (0,-4) {\vdots};
    \node[main node, minimum size=1cm] (bm) at (0,-5) {\small $n_b$ };
    
    \foreach \x in {1,i,m}
        \draw (b\x) -- ++(-2,0);
    
    \node[red, left] at (-2,-1) {$y_1 = u_1(x_1)$};
    \node[blue, left] at (-2,-3) {$y_i = u_i(x_i)$};
    \node[teal, left] at (-2,-5) {$y_m = u_{n_b}(x_{n_b})$};

    
    \node (goods) at (4,0) {goods};
    \node[main node, minimum size=1cm] (g1) at (4,-1) {\tiny $n_b + 1$ };
    \node at (4,-2) {\vdots};
    \node[main node, minimum size=1cm] (gj) at (4,-3) {\small $j$ };
    \node at (4,-4) {\vdots};
    \node[main node, minimum size=1cm] (gn) at (4,-5) {\small $n_g$ };
    
    \foreach \x in {1,j,n}
    {
        \draw[red] (g\x) -- (b1);
        \draw[blue] (g\x) -- (bi);
        \draw[teal] (g\x) -- (bm);
    }
    
    
    \draw[->] (6,-1) node[right] {$y_{n_b+1} \ge -1$} -- (g1);
    \draw[->] (6,-3) node[right] {$y_{j} \ge -1$} -- (gj);
    \draw[->] (6,-5) node[right] {$y_{n_b+n_g} \ge -1$} -- (gn);
    
a    \end{tikzpicture}
    \caption{Network representation of the linear Fisher market model, where
    we aim to maximize the utility of the net flows $y$. Colored edges represent 
    the $n_b+1$ hyperedges connecting each buyer to all the goods, so each of 
    these edges is incident to $n_g + 1$ vertices. Flows on the right indicate
    that there is at most one unit of each good to be divided among the buyers.}
    \label{fig:fisher-mkt}
\end{figure}

\paragraph{The convex flow problem.}
This problem can be turned into a convex flow problem in a number of ways.
Here, we follow V\'{e}gh~\cite{vegh2014concave} and use a similar construction.
First, we represent both the $n_b$ buyers and $n_g$ goods as nodes on the graph: 
the buyers are labeled as nodes $1, \dots, n_b$, and the goods are 
labeled as nodes $n_b+1, \dots, n_b + n_g$. For each buyer $i = 1, \dots, n_b$, 
we introduce a hyperedge connecting buyer $i$ to all goods $j = n_b+1, \dots, n_b + n_g$.
We denote the flow along this edge by $x_i \in \reals^{n_g+1}$, where $(x_i)_j$ 
is the amount of the $j$th good bought by the $i$th buyer, and $(x_i)_{n_g+1}$ 
denotes the amount of utility that buyer $i$ receives from this basket of goods 
denoted by the first $n_g$ entries of $x_i$. The flows on this edge are given by 
the convex set
\begin{equation}
    \label{eq:app-market-equilibrium-Ti}
    T_i = \left\{ 
        (z, t) \in \reals^{n_g} \times \reals 
        \mid 
        -\ones \le z \le 0,~~
        t \le u_i(-z)
    \right\}.
\end{equation}
This set converts the `goods' flow into a utility flow at each buyer node $i$.
This setup is depicted in figure~\ref{fig:fisher-mkt}. Since the indices $1,
\dots, n_b$ of the net flow vector $y \in \reals^{n_g + n_b}$ correspond to
the buyers and the elements $n_b+1, \dots, n_g + n_b$ to correspond to the
goods, the utility function above can be written
\begin{equation}
    \label{eq:app-market-equilibrium-U}
    U(y) = \sum_{i=1}^{n_b} b_i\log(y_i) - I(y_{n_b+1:n_b+n_g} \ge -\ones),
\end{equation}
which includes the implicit constraint that at most one unit of each good can
be sold in the indicator function $I$, above, and where $y_{n_b+1:n_b+n_g}$ is
a vector containing only the $(n_b+1)$st to the $(n_b + n_g)$th entries of $y$.
Since $U$ is nondecreasing and concave in $y$, $V_i = 0$, and the $T_i$ are
convex, we know that~\eqref{eq:app-market-equilibrium} is a special case
of the convex flow problem~\eqref{eq:main-problem}.

\paragraph{Related problems.}
A number of other resource allocation can be cast as generalized network flow
problems. For example, Agrawal et al.~\cite{agrawal2022allocation} consider a
price adjustment algorithm for allocating compute resources to a set of jobs, 
and Schutz et al.~\cite{schutz2009supply} do the same for supply chain problems.
In many problems, network structure implicitly appears if we are forced to make
decisions over time or over decision variables which directly interact only with 
a small subset of other variables.

\subsection{Routing orders through financial exchanges}\label{sec:app-cfmm}
Financial asset networks are also well-modeled by convex network flows. If each
asset is a node and each market between assets is an edge between the corresponding nodes,
we expect the edge input-output functions to be concave, as the price of an
asset is nondecreasing in the quantity purchased. In many markets, this
input-output function is probabilistic; the state of the market when the order
`hits' is unknown due to factors such as information latency, stale orders, and
front-running. However, in certain batched exchanges, including
\emph{decentralized exchanges} running on public blockchains, this state can be
known in advance. We explore the order routing problem in this decentralized
exchange setting.

\paragraph{Decentralized exchanges and automated market makers.}
Automated market makers have reached mass adoption after being implemented as
decentralized exchanges on public blockchains. These exchanges (including Curve
Finance~\cite{egorov2019stableswap}, Uniswap~\cite{adams2020uniswap}, and
Balancer~\cite{balancer}, among others) have facilitated trillions of
dollars in cumulative trading volume since 2019 and maintain a collective daily 
trading volume of several billion dollars. These exchanges are almost all 
implemented as constant function market markers (CFMMs)~\cite{angerisImprovedPriceOracles2020,angeris2023geometry}. 
In CFMMs, liquidity providers contribute reserves of assets. Users can then 
trade against these reserves by tendering a basket of assets in exchange for
another basket. CFMMs use a simple rule for accepting trades: a trade is only 
valid if the value of a given function at the post-trade reserves is equal to 
the value at the pre-trade reserves. This function is called the trading 
function and gives CFMMs their name.

\paragraph{Constant function market makers.}
A CFMM which allows $r$ assets to be traded is defined by two properties: its
reserves $R \in \reals^r$, which denotes the amount of each asset available to
the CFMM, and its trading function $\phi : \reals^r \to \reals$, which
specifies its behavior and includes a fee parameter $0 < \gamma \le 1$, where
$\gamma = 1$ denotes no fee. We assume that $\phi$ is concave and
nondecreasing. Any user is allowed to submit a trade to a CFMM, which we write
as a vector $z \in \reals^r$, where positive entries denote values to be
received from the CFMM and negative entries denote values to be tendered to the
CFMM. (For example, if $r=2$, then a trade $z = (-1, 10)$ would denote that the user
wishes to tender $1$ unit of asset $1$ and receive $10$ units of asset
$2$.) The submitted trade is then accepted if the following condition holds:
\[
    \phi(R - \gamma z_- - z_+) \ge \phi(R),
\]
and $R - \gamma z_- - z_+ \ge 0$. Here, we denote $z_+$ to be the `elementwise
positive part' of $x$, \ie, $(z_+)_j = \max\{z_j, 0\}$ and $z_-$ to be the
`elementwise negative part' of $x$, \ie, $(z_-)_j = \min\{z_j, 0\}$ for every
asset $j=1, \dots, r$. Note that, since $\phi$ is concave, the set of acceptable
trades is a convex set:
\[
    T = \{z \in \reals^r \mid \phi(R - \gamma z_- - z_+) \ge \phi(R)\},
\]
as we can equivalently write it as
\[
    T = \{z \in \reals^r \mid \phi(R + \gamma u - v) \ge \phi(R), ~ u, v \ge 0, ~ z = v - u \},
\]
which is easily seen to be a convex set since $\phi$ is a concave function.

\paragraph{Net trade.}
Consider a collection of $m$ CFMMs, each of which trades a subset of $n$
possible assets. Denoting the trade with the $i$th CFMM by $x_i$, which must lie
in the convex set $T_i$, we can write the net trade across all markets by $y \in
\reals^n$, where
\[
    y = \sum_{i=1}^m A_i x_i, \quad \text{and} \quad x_i \in T_i.
\]
If $y_j > 0$, we receive some amount of asset $j$ after executing all trades
$\left\{x_i\right\}$. On the other hand, if $y_j < 0$, we tender some of asset
$j$ to the network.

\paragraph{Optimal routing.}
Finally, we denote the trader's utility of the network trade vector by $U: 
\reals^n \to \reals \cup \{-\infty\}$, where infinite values encode constraints. 
We assume that this function is concave and nondecreasing. We can choose $U$ to
encode several important actions in markets, including liquidating a portfolio, 
purchasing a basket of assets, and finding arbitrage. For example, if we wish to
find risk-free arbitrage, we may take
\[
    U(y) = c^Ty - I(y \ge 0),
\]
for some vector of prices $c \in \reals^n$. See~\cite[\S
5.2]{angerisConstantFunctionMarket2021} for several additional examples.
Letting $V_i = 0$ for all $i = 1, \dots, m$, it's clear that the optimal
routing problem in CFMMs is a special case of~\eqref{eq:main-problem}.

\section{The dual problem and flow prices}\label{sec:dual-problem}
The remainder of the paper focuses on efficiently solving the convex flow
problem~\eqref{eq:main-problem}. This problem has only one constraint coupling
the edge flows and the net flow variables. As a result, we turn to dual
decomposition methods~\cite{boyd2007notes,bertsekasNonlinearProgramming2016}.
The general idea of dual decomposition methods is to solve the original problem
by splitting it into a number of subproblems that can be solved quickly and
independently. In this section, we will design a decomposition method that
parallelizes over all edges and takes advantage of structure present in the
original problem. This decomposition allows us to quickly evaluate the dual
function and a subgradient. Importantly, our decomposition method also provides
a clear programmatic interface to specify and solve the convex flow problem.

\subsection{Dual decomposition}
To get a dual problem, we introduce a set of (redundant) additional variables
for each edge and rewrite~\eqref{eq:main-problem} as
\[
    \begin{aligned}
        & \text{maximize} && {\textstyle U(y) + \sum_{i=1}^m V_i(x_i)}\\
        & \text{subject to} && {\textstyle y = \sum_{i=1}^m A_i x_i}\\
        &&& x_i = \tilde x_i, ~~ \tilde x_i \in T_i, \quad i=1, \dots, m,
    \end{aligned}
\]
where we added the `dummy' variables $\tilde x_i \in \reals^{n_i}$ for $i=1,
\dots, m$. Next, we pull the constraint $\tilde x_i \in T_i$ for $i=1, \dots,
m$ into the objective by defining the indicator function
\[
    I_i(\tilde x_i) = \begin{cases}
        0 & \tilde x_i \in T_i\\
        +\infty & \text{otherwise}.
    \end{cases}
\]
This rewriting gives the augmented problem,
\begin{equation}
    \label{eq:augmented-problem}
    \begin{aligned}
        & \text{maximize} && {\textstyle U(y) + \sum_{i=1}^m \left(V_i(x_i) - I_i(\tilde x_i)\right)}\\
        & \text{subject to} && {\textstyle y = \sum_{i=1}^m A_i x_i}\\
        &&& x_i = \tilde x_i, \quad i=1, \dots, m,
    \end{aligned}
\end{equation}
with variables $x_i, \tilde x_i \in \reals^{n_i}$ for $i = 1, \dots, m$
and $y \in \reals^n$.
The Lagrangian~\cite[\S 5.1.1]{cvxbook} of this problem is then
\begin{equation}\label{eq:lagrangian}
L(y, x, \tilde x, \nu, \eta) =
    U(y) - \nu^Ty + \sum_{i=1}^m \left(V_i(x_i) + (A_i^T\nu - \eta_i)^Tx_i\right) + \sum_{i=1}^m \left(\eta_i ^T\tilde x_i - I_i(\tilde x_i) \right),
\end{equation}
where we have introduced the dual variables $\nu \in \reals^n$ for the net flow
constraint and $\eta_i \in \reals^{n_i}$ for $i=1, \dots, m$ for each of the
individual edge constraints in~\eqref{eq:augmented-problem}. (We will write $x$, $\tilde x$,
and $\eta$ as shorthand for $\{x_i\}$, $\{\tilde x_i\}$, and $\{\eta_i\}$, respectively.)
It's easy to see that the Lagrangian is separable over the primal variables $y$,
$x$, and $\tilde x$.

\paragraph{Dual function.}
Maximizing the Lagrangian~\eqref{eq:lagrangian} over the primal variables $y$,
$x$, and $\tilde x$ gives the dual function
\[
    \begin{aligned}
        g(\nu, \eta) 
        &= \sup_{y} (U(y) - \nu^Ty) + \sum_{i=1}^m \left(
            \sup_{x} (V_i(x) - (A_i^T\nu - \eta_i)^Tx) 
            + \sup_{\tilde x_i \in T_i} \eta_i^T\tilde x_i
        \right).
    \end{aligned}
\]
To evaluate the dual function, we must solve three subproblems, each
parameterized by the dual variables $\nu$ and $\eta$. We denote the optimal
values of these problems, which depend on the $\nu$ and $\eta$, by $\bar U$,
$\bar V_i$, and $f_i$:
\begin{subequations}
    \label{eq:subproblems}
    \begin{align}
        \label{eq:subproblem-U}
        \bar U(\nu) &= \sup_{y} (U(y) - \nu^Ty), \\
        \label{eq:subproblem-V}
        \bar V_i(\xi) &= \sup_{x_i} (V_i(x_i) - \xi^Tx_i), \\
        \label{eq:subproblem-f}
        f_i(\tilde \eta) &= \sup_{\tilde x_i \in T_i} \tilde \eta_i^T \tilde x_i.
    \end{align}
\end{subequations}
The functions $\bar U$ and $\{\bar V_i\}$ are essentially the Fenchel
conjugate~\cite[\S3.3]{cvxbook} of the corresponding $U$ and $\{V_i\}$.
Closed-form expressions for $\bar U$ and the $\{\bar V_i\}$ are known for many
practical functions $U$ and $\{V_i\}$. Similarly, the functions $\{f_i\}$ are
the support functions~\cite[\S13]{rockafellar1970convex} for the sets $T_i$. For future
reference, note that the $\bar U$, $\bar V_i$, and $f_i$ are convex, as they
are the pointwise supremum of a family of affine functions, and may take on
value $+\infty$, which we interpret as an implicit constraint. We can rewrite
the dual function in terms of these functions~\eqref{eq:subproblems} as
\begin{equation}
    \label{eq:dual-function}
    g(\nu, \eta) = \bar U(\nu) + \sum_{i=1}^m \left(\bar V_i(\eta_i - A_i^T\nu) + f_i(\eta_i)\right).
\end{equation}
Our ability to quickly evaluate these functions and their gradients governs the
speed of any optimization algorithm we use to solve the dual problem.
The dual function~\eqref{eq:dual-function} also has very clear structure: the `global'
dual variables $\nu$ are connected to the `local' dual variables $\eta$, only
through the functions $\{\bar V_i\}$. If the $\bar V_i$ were all affine 
functions, then the problem would be separable over $\nu$ and each $\eta_i$.

\paragraph{Dual variables as prices.}
Subproblem~\eqref{eq:subproblem-U} for evaluating $\bar U(\nu)$ has a simple
interpretation: if the net flows $y$ have per-unit prices $\nu \in \reals^n$,
find the maximum net utility, after removing costs, over all net flows. (There
need not be feasible edge flows $x$ which correspond to this net flow.)
Assuming $U$ is differentiable, a $y$ achieving this maximum satisfies
\[
    \nabla U(y) = \nu,
\]
\ie, the marginal utilities for flows $y$ are equal to the prices $\nu$. (A
similar statement for a non-differentiable $U$ follows directly from
subgradient calculus.) The subproblems over the $V_i$~\eqref{eq:subproblem-V}
have a similar interpretation as utility maximization problems.

On the other hand, subproblem~\eqref{eq:subproblem-f} for evaluating
$f_i(\tilde\eta)$ can be interpreted as finding a most `valuable' allowable flow over
edge $i$. In other words, if there exists an external, infinitely liquid
reference market where we can buy or sell flows $x_i$ for prices $\tilde \eta \in
\reals^{n_i}$, then $f_i(\tilde \eta)$ gives the highest net value of any allowable
flow $x_i \in T_i$. Due to this interpretation, we will refer
to~\eqref{eq:subproblem-f} as the \emph{arbitrage problem}. This price
interpretation is also a natural consequence of the optimality conditions for
this subproblem. The optimal flow $x^0$ is a point in $T_i$ such that there
exists a supporting hyperplane to $T_i$ at $x^0$ with slope $\tilde \eta$. In other
words, for any small deviation $\delta \in \reals^{n_i}$, if $x^0 + \delta \in
T_i$, then 
\[
    \tilde\eta^T(x^0 + \delta) \le \tilde\eta^Tx^0 \implies \tilde\eta^T\delta \le 0.
\]
If, for example, $\delta_j$ and $\delta_k$ are the only two nonzero entries of
$\delta$, we would have
\[
    \delta_j \le -\frac{\tilde\eta_k}{\tilde\eta_j}\delta_k,
\]
so the exchange rate between $j$ and $k$ is at most $\tilde\eta_j/\tilde\eta_k$.
This observation lets us interpret the dual variables $\eta$ as `marginal prices'
on each edge, up to a constant multiple. With this interpretation, we will soon 
see that the function $\bar V_i$ also connects the `local prices' $\eta_i$ on 
edge $i$ to the `global prices' $\nu$ over the whole network.

\paragraph{Duality.}
An important consequence of the definition of the dual function is weak
duality~\cite[\S5.2.2]{cvxbook}. Letting $p^\star$ be an optimal value for the
convex flow problem~\eqref{eq:main-problem}, we have that
\begin{equation}\label{eq:weak-duality}
    g(\nu, \eta) \ge p^\star.
\end{equation}
for every possible choice of $\nu$ and $\eta$. 
An important (but standard) result in convex optimization states that there
exists a set of prices $(\nu^\star, \eta^\star)$ which actually achieve
the bound:
\[
    g(\nu^\star, \eta^\star) = p^\star,
\]
under mild conditions on the problem data~\cite[\S5.2]{cvxbook}. One such
condition is if all the $T_i$'s are affine sets, as
in~\S\ref{sec:app-max-flow}. Another is Slater's condition: if there exists a
point in the relative interior of the feasible set, \ie, if the set
\[
    \left(\sum_{i=1}^m \relint\left(A_i \left(T_i \cap \dom V_i\right)\right)\right) \cap \relint\dom U
\]
is nonempty. (We have used the fact that the $A_i$ are one-to-one projections.)
These conditions are relatively technical but almost always hold in practice.
We assume they hold for the remainder of this section.

\subsection{The dual problem}
The dual problem is then to find a set of prices $\nu^\star$ and $\eta^\star$
which saturate the bound~\eqref{eq:weak-duality} at equality; or, equivalently, 
the problem is to find a set of prices that minimize the dual function $g$. 
Using the definition of $g$ in~\eqref{eq:dual-function}, we may write this 
problem as
\begin{equation}
    \label{eq:dual-problem-impl}
    \begin{aligned}
        & \text{minimize} && \bar U(\nu) + \sum_{i=1}^m \left(\bar V_i(\eta_i - A_i^T\nu) + f_i(\eta_i)\right),
    \end{aligned}
\end{equation}
over variables $\nu$ and $\eta$. The dual problem is a convex optimization
problem since $\bar U$, $\bar V_i$, and $f_i$ are all convex functions. For
fixed $\nu$, the dual problem~\eqref{eq:dual-problem} is also separable over
the dual variables $\eta_i$ for $i=1, \dots, m$; we will later use this fact to
speed up solving the problem by parallelizing our evaluations of each $\bar V_i$
and $f_i$.

\paragraph{Implicit constraints.} The `unconstrained'
problem~\eqref{eq:dual-problem-impl} has implicit constraints due to the fact that
the $U$ and $V_i$ are nondecreasing functions. More specifically, if $U$
is nondecreasing and $U(0) < \infty$, then, if $\nu_i < 0$, we have
\[
    \bar U(\nu) = \sup_y (U(y) - \nu^Ty) \ge U(te_i) - t\nu_i \ge U(0) - t\nu_i \to \infty,
\]
as $t \uparrow \infty$. Here, in the first inequality, we have chosen $y =
te_i$, where $e_i$ is the $i$th unit basis vector. This implies that $\bar
U(\nu) = \infty$ if $\nu \not \ge 0$, which means that $\nu \ge 0$ is an
implicit constraint. A similar proof shows that $\bar V_i(\xi) = \infty$
if $\xi \not\ge 0$. Adding both implicit constraints as explicit
constraints gives the following constrained optimization problem:
\begin{equation}\label{eq:dual-problem}
    \begin{aligned}
        & \text{minimize} && \bar U(\nu) + \sum_{i=1}^m \left(\bar V_i(\eta_i - A_i^T\nu) + f_i(\eta_i)\right)\\
        & \text{subject to} && \nu \ge 0, \quad \eta_i \ge A_i^T\nu, \quad i=1, \dots, m.
    \end{aligned}
\end{equation}
Note that this implicit constraint exists even if $U(0) = \infty$; we only
require that the domain of $U$ is nonempty, \ie, that there exists some $y$
with $U(y) < \infty$, and similarly for the $V_i$. The result follows from a
nearly-identical proof. This fact has a simple interpretation in the context
of utility maximization as discussed previously: if we have a nondecreasing
utility function and are paid to receive some flow, we will always choose to 
receive more of it.

In general, the rewriting of problem~\eqref{eq:dual-problem-impl} into
problem~\eqref{eq:dual-problem} is useful since, in practice, $\bar U(\nu)$ is
finite (and often differentiable) whenever $\nu \ge 0$; a similar thing is true
for the functions $\{V_i\}$. Of course, this need not always be true, in which case the
additional implicit constraints need to be made explicit in order to use
standard, off-the-shelf solvers for types of problems.

\paragraph{Optimality conditions.}
Let $(\nu^\star, \eta^\star)$ be an optimal point for the dual problem, and
assume that $g$ is differentiable at this point.
The optimality conditions for the dual problem are then
\[
   \nabla g(\nu^\star, \eta^\star) = 0.
\]
(The function $g$ need not be differentiable, in which case a similar argument 
holds using subgradient calculus.) For a differentiable $\bar U$, we
have that
\[
   \nabla_\nu \bar U(\nu^\star) = -y^\star  
\] 
where $y^\star$ is the optimal point for subproblem~\eqref{eq:subproblem-U}.
For a differentiable $\bar V_i$, we have that
\[
   \begin{aligned}
       \nabla_\nu \bar V_i(\eta_i^\star - A_i^T\nu^\star) = A_ix_i^\star, \\
       \nabla_{\eta_i} \bar V_i(\eta_i^\star - A_i^T\nu^\star) = -x_i^\star.
   \end{aligned}
\]
where $x_i^\star$ is the optimal point for subproblem~\eqref{eq:subproblem-V}.
Finally, we have that 
\[
   \nabla f_i(\eta_i^\star) = \tilde x_i^\star,
\]
where $\tilde x_i^\star$ is the optimal point for
subproblem~\eqref{eq:subproblem-f}. Putting these together with the definition
of the dual function~\eqref{eq:dual-function}, we recover primal feasibility at
optimality:
\begin{equation}
   \label{eq:gradient}
   \begin{aligned}
       \nabla_\nu g(\nu^\star, \eta^\star)
       &= y^\star - \sum_{i=1}^m A_i x_i^\star = 0, \\
       \nabla_{\eta_i} g(\nu^\star, \eta^\star)
       &= x_i^\star - \tilde x_i^\star = 0, \quad i=1, \dots, m.
   \end{aligned}
\end{equation}
In other words, by choosing the `correct' prices $\nu^\star$ and $\eta^\star$
(\ie, those which minimize the dual function), we find that the optimal
solutions to the subproblems in~\eqref{eq:subproblems} satisfy the resulting
coupling constraints, when the functions in~\eqref{eq:subproblems} are all
differentiable at the optimal prices $\nu^\star$ and $\eta^\star$. This, in
turn, implies that the $\{x_i^\star\}$ and $y^\star$ are a solution to the
original problem~\eqref{eq:main-problem}. In the case that the functions are
not differentiable, there might be many optimal solutions to the subproblems
of~\eqref{eq:subproblems}, and we are only guaranteed that at least one of
these solutions satisfies primal feasibility. We give some heuristics to handle
this latter case in~\S\ref{sec:dual-primal-feas}.

\paragraph{Dual optimality conditions.}
For problem~\eqref{eq:dual-problem}, if the functions $U$ and
$V_i$ are differentiable at optimality, the dual conditions state that
\begin{equation}\label{eq:dual-conditions}
   \begin{aligned}
       \nabla U(y^\star) 
       &= \nu^\star, \\ 
       \nabla V_i(x_i^\star) 
       &= \eta_i^\star - A_i^T\nu^\star, \quad i = 1, \dots, m \\
       \eta_i^\star 
       &\in \mathcal{N}_i(\tilde x_i^\star), \quad i = 1, \dots, m,
   \end{aligned}
\end{equation}
where $\mathcal{N}_i(x)$ is the normal cone for set $T_i$ at the point $x$,
defined as
\[
    \mathcal{N}_i(x) = \{u \in \reals^{n_i} \mid u^Tx \ge u^Tz ~~ \text{for all} ~ z \in T_i\}.
\]
Note that, since $U$ is nondecreasing, then, if $U$ is differentiable, its
gradient must be nonnegative, which includes the implicit constraints
in~\eqref{eq:dual-problem}. (A similar thing is true for the $V_i$.) A similar
argument holds in the case that $U$ and the $V_i$ are not differentiable at
optimality, via subgradient calculus, and the implicit constraints are
similarly present.

\paragraph{Interpretations.} 
The dual optimality conditions~\eqref{eq:dual-conditions} each have lovely
economic interpretations. In particular, they state that, at optimality, the
marginal utilities from the net flows $y^\star$ are equal to the prices
$\nu^\star$, the marginal utilities from the edge flows $x_i^\star$ are equal to
the difference between the local prices $\eta_i^\star$ and the global prices
$\nu^\star$, and the local prices $\eta_i^\star$ are a supporting hyperplane of
the set of allowable flows $T_i$. This interpretation is natural in problems
involving markets for goods, such as those discussed in~\S\ref{sec:app-opf}
and~\S\ref{sec:app-cfmm}, where one may interpret the normal cone
$\mathcal{N}_i(x_i^\star)$ as the no-arbitrage cone for the market $T_i$: any
change in the local prices $\eta_i^\star$ such that $\eta_i^\star$ is still in
the normal cone $\mathcal{N}_i(\tilde x_i^\star)$ does not affect our action
$\tilde x^\star_i$ with respect to market $i$.

We can also interpret the dual problem~\eqref{eq:dual-problem} similarly to the
primal. Here, each subsystem has local `prices' $\eta_i$ and is described by
the functions $f_i$ and $\bar V_i$, which implicitly include the edge flows and
associated constraints. The global prices $\nu$ are associated with a cost
function $\bar U$, which implicitly depends on (potentially infeasible) net
flows, $y$. The function $\bar V_i$ may be interpreted as a convex cost
function of the difference between the local subsystem prices $\eta_i$ and the
global prices $\nu$. At optimality, this difference will be equal to the
marginal utility gained from the optimal flows in subsystem $i$, and the global
prices will be equal to the marginal utility of the overall system.

Finally, we note that the convex flow problem~\eqref{eq:main-problem} and its 
dual~\eqref{eq:dual-problem-impl} examined in this section are closely related 
to the extended monotropic programming problem~\cite{bertsekas2008extended}. We 
make this connection explicit in appendix~\ref{app:extended-monotropic}. 
Importantly, the extended monotropic programming problem is self-dual, whereas
the convex flow problem is not evidently self-dual. This observation suggests an
interesting avenue for future work.

\subsection{Zero edge utilities}\label{sec:special-cases}
An important special case is when the edge flow utilities are zero, \ie, if $V_i = 0$ for $i = 1, \dots, m$.
In this case, the convex flow problem reduces to the routing problem
discussed in~\S\ref{sec:app-cfmm}, originally presented
in~\cite{angeris2022optimal,diamandis2023efficient} in the context of constant
function market makers~\cite{angerisImprovedPriceOracles2020}. Note that $\bar
V_i$ becomes
\[
    \bar V_i(\xi_i) = \begin{cases}
        0 & \xi_i = 0\\
        +\infty &\text{otherwise},
    \end{cases}
\]
which means that the dual problem is
\[
    \begin{aligned}
        & \text{minimize} && {\textstyle \bar U(\nu) + \sum_{i=1}^m f_i(\eta_i)}\\
        & \text{subject to} && \eta_i = A_i^T\nu, \quad i=1, \dots, m.
    \end{aligned}
\]
This equality constraint can be interpreted as ensuring that the local prices
for each node are equal to the global prices over the net flows of the network.
If we substitute $\eta_i = A_i^T\nu$ in the objective, we have
\begin{equation}\label{eq:zero-edge-problem}
    \begin{aligned}
        & \text{minimize} && \bar U(\nu) + \sum_{i=1}^m f_i(A_i^T\nu),
    \end{aligned}
\end{equation}
which is exactly the dual of the optimal routing problem, originally presented
in~\cite{diamandis2023efficient}. In the case of constant function market
makers (see \S\ref{sec:app-cfmm}), we interpret the subproblem of computing the
value of $f_i$, at some prices $A_i^T\nu$, as finding the optimal arbitrage with
the market described by $T_i$, given `true' (global) asset prices $\nu$. 

\paragraph{Problem size.} Because this problem has only $\nu$ as a variable,
which is of length $n$, this problem is often much smaller than the original
dual problem of~\eqref{eq:dual-problem}. Indeed, the number of variables in the
original dual problem is $n + \sum_{i=1}^m n_i \ge n + 2m$, whereas this
problem has exactly $n$ variables. (Here, we have assumed that the feasible
flow sets lie in a space of at least two dimensions, $n_i \ge 2$.) This special
case is very common in practice and identifying it often leads to significantly
faster solution times, as the number of edges in many practical networks is much 
larger than the total number of nodes, \ie, $m \gg n$.

\paragraph{Example.} Using this special case, is easy to show that the dual for
the maximum flow problem~\eqref{eq:max-flow}, introduced
in~\S\ref{sec:app-max-flow}, is the minimum cut problem, as expected from the
celebrated result of~\cite{dantzig1955max,elias1956note,ford1956maximal}. 
Recall from~\S\ref{sec:app-max-flow} that
\[
    U(y) = y_n - I_S(y), \qquad T_i = \{z \in \reals^2 \mid 0 \le z_2 \le b_i, ~~ z_1+z_2=0\},
\]
where $S = \{y \in \reals^n \mid y_1+y_n \ge 0, ~ y_i \ge 0 ~ \text{for all} ~
i \ne 1, n\}$ and $b_i \ge 0$ is the maximum allowable flow across edge $i$.
Using the definitions of $\bar U$ and $f_i$ in~\eqref{eq:subproblems}, we can 
easily compute $\bar U(\nu)$,
\[
    \bar U(\nu) = \sup_{y \in S} \left( y_n - \nu^T y\right) = \begin{cases}
        0 & \nu_n \ge 1, ~ \nu_n  - \nu_1 = 1, ~ \nu_i \ge 0 ~~ \text{for all} ~ i\ne 1, n\\
        +\infty & \text{otherwise},
    \end{cases}
\]
and $f_i(\eta)$,
\[
    f_i(\eta) = \sup_{z \in T_i} \eta^Tz = b_i(\eta_2 - \eta_1)_+,
\]
where we write $(w)_+ = \max\{w, 0\}$. Using the special case of the problem
when we have zero edge utilities~\eqref{eq:zero-edge-problem} and adding the
constraints gives the dual problem 
\[
\begin{aligned}
& \text{minimize} && \sum_{i=1}^m b_i((A_i^T\nu)_2 - (A_i^T\nu)_1)_+ \\
& \text{subject to} && \nu_n - \nu_1 = 1 \\
&&& \nu_n \ge 1, ~\nu_{i} \ge 0, ~\text{for all} ~ i\ne 1, n.
\end{aligned}
\]
Note that this problem is $1$-translation invariant: the problem has the same
objective value and remains feasible if we replace any feasible $\nu$ by $\nu +
\alpha \ones$ for any constant $\alpha$ such that $\nu_n +\alpha \ge 1$. Thus,
without loss of generality, we may always set $\nu_1 = 0$ and $\nu_n = 1$. We
then use an epigraph transformation and introduce new variables for each edge,
$t \in \reals^m$, so the problem becomes
\[
\begin{aligned}
& \text{minimize} && b^Tt \\
& \text{subject to} && (A_i^T\nu)_1 - (A_i^T\nu)_2 \le t_i, \quad i = 1, \dots n \\
&&& \nu_n = 1, ~ \nu_1 = 0\\
&&& t\ge 0, ~ \nu \ge 0.
\end{aligned}
\]
The substitution of $\nu_n = 1$ and $\nu_1 = 0$ in the first constraint 
recovers a standard formulation of the minimum cut problem (see, \eg,~\cite[\S3]{dantzig1955max}).

\section{Solving the dual problem.}\label{sec:dual}\label{sec:algorithm}
The dual problem~\eqref{eq:dual-problem} is a convex optimization problem that 
is easily solvable in practice, even for very large $n$ and $m$. For small problem
sizes, we can use an off-the-self solver, such as such as SCS~\cite{odonoghueConicOptimizationOperator2016},
Hypatia~\cite{coey2021solving}, or Mosek~\cite{mosek}, to tackle the convex 
flow problem~\eqref{eq:main-problem} directly; however, these methods, which
rely on conic reformulations, destroy problem structure and may be unacceptably
slow for large problem sizes. The dual problem, on the other hand, preserves 
this structure, so our approach is to solve this dual problem. 

\paragraph{A simple transformation.}
For the sake of exposition, we will introduce the new variable $\mu = (\nu, \eta)$
and write the dual problem~\eqref{eq:dual-problem} as 
\[
\begin{aligned}
    & \text{minimize} && g(\mu) \\
    & \text{subject to} && F \mu \ge 0,
\end{aligned}
\]
where $F$ is the constraint matrix
\[
    F = \bmat{
        I & 0 & \cdots & 0 \\
        -A_1^T & I & \cdots & 0 \\
        \vdots & \vdots & \ddots & \vdots \\
        -A_m^T & 0 & \cdots & I
    }.
\]
Since the matrix $F$ is lower triangular with a diagonal that has no nonzero
entries, the matrix $F$ is invertible. Its inverse is given by
\[
    F^{-1} = \bmat{
        I & 0 & \cdots & 0 \\
        A_1^T & I & \cdots & 0 \\
        \vdots & \vdots & \ddots & \vdots \\
        A_m^T & 0 & \cdots & I
    },
\]
which can be very efficiently applied to a vector. With these matrices defined,
we can rewrite the dual problem as
\begin{equation}
    \label{eq:dual-problem-nonnegative-cone}
    \begin{aligned}
        & \text{minimize} && g(F^{-1}\tilde\mu) \\
        & \text{subject to} && \tilde\mu \ge 0,
    \end{aligned}
\end{equation}
where $\tilde \mu = F\mu$. Note that the matrix $F$ preserves the separable
structure of the problem: edges are only directly coupled with their adjacent 
nodes.

\paragraph{An efficient algorithm.}
After this transformation, we can apply any first-order method that handles
nonnegativity constraints to solve~\eqref{eq:dual-problem-nonnegative-cone}. We
use the quasi-Newton algorithm
L-BFGS-B~\cite{byrd1995limited,zhu1997algorithm,morales2011remark}, which has
had much success in practice. This algorithm only requires the ability to
evaluate the function $g$ and its gradient $\nabla g$ at a given point $\mu =
F^{-1}\tilde \mu$. The function $g(\mu)$ and its gradient $\nabla g(\mu)$ are, 
respectively, a sum of the optimal values and a sum of the optimal points for 
the subproblems~\eqref{eq:subproblems}.

\paragraph{Interface.}
The use of a first-order method suggests a natural interface for specifying the
convex flow problem~\eqref{eq:main-problem}. By definition, the function $g$ is
easy to evaluate if the subproblems~\eqref{eq:subproblem-U},
\eqref{eq:subproblem-V}, and \eqref{eq:subproblem-f} are easy to evaluate.
Given a way to evaluate the functions $\bar U$, $\bar V_i$, and $f_i$, and to
get the values achieving the suprema in these subproblems, we can easily
evaluate the function $g$ via~\eqref{eq:dual-function} and its gradient $\nabla
g$ via~\eqref{eq:gradient}, which we write below:
\[
    \begin{aligned}
        \nabla_\nu g(\nu, \eta) &= y^\star - \sum_{i=1}^m A_i x_i^\star \\
        \nabla_{\eta_i} g(\nu, \eta) &= x_i^\star - \tilde x_i^\star, \quad i=1, \dots, m.
    \end{aligned}
\]
(Here, as before, $y^\star$ and the $\{x_i^\star\}$ and $\{\tilde x_i^\star\}$
are the optimal points for the subproblems~\eqref{eq:subproblems}, evaluated at
$\eta$ and $\nu$.) Often $\bar U$ and $\bar V_i$, which are closely related to
conjugate functions, have a closed form expression. In general, evaluating the
support function $f_i$ requires solving a convex optimization problem. However,
in many practical scenarios, this function either has a closed form expression
or there is an easy and efficient iterative method for evaluating it. We discuss a method
for quickly evaluating this function in the special case of two-node edges
in~\S\ref{sec:dual-two-node} and implement more general subproblem solvers for
the examples in~\S\ref{sec:numerical-examples}.

\paragraph{Parallelization.}
The evaluation of $g(\nu, \eta)$ and $\nabla g(\nu, \eta)$ can be parallelized
over all edges $i = 1, \dots, m$. The bulk of the computation is in evaluating
the arbitrage subproblem $f_i$ for each edge $i$. The optimality conditions of
the subproblem suggest a useful shortcut: if the vector $0$ is in the normal
cone of $T_i$ at $\tilde x_i^\star$, then the zero vector is a solution to the
subproblem, \ie, the edge is not used (see~\eqref{eq:dual-conditions}). Often,
this condition is not only easier to evaluate than the subproblem itself, but
also has a nice interpretation. For example, in the case of financial markets,
this condition is equivalent to the prices $\eta$ being within the bid-ask
spread of the market. We will see that, in many examples, these subproblems are
quick to solve and have further structure that can be exploited.

\subsection{Two-node edges}\label{sec:dual-two-node}
In many applications, edges are often between only two vertices. Since these
edges are so common, we will discuss how the additional structure allows the
arbitrage problem~\eqref{eq:subproblem-f} to be solved quickly for such edges.
Some practical examples will be given later in the numerical examples
in~\S\ref{sec:numerical-examples}. In this section, we will drop the index $i$
with the understanding that we are referring to the flow along a particular
edge. Note that much of this section is similar to some of the authors' previous
work~\cite[\S3]{diamandis2023efficient}, where two-node edges were explored in
the context of the CFMM routing problem, also discussed in~\S\ref{sec:app-cfmm}.

\subsubsection{Gain functions}
To efficiently deal with two-node edges, we will consider their \emph{gain 
functions}, which denote the maximum amount of output one can receive given a
specified input. Note that our gain function is equivalent to the concave gain 
functions introduced in~\cite{shigeno2006maximum}, and, in the case of asset 
markets, the forward exchange function introduced in~\cite{angerisConstantFunctionMarket2021}.
In what follows, the gain function $h: \reals \to \reals \cup \{-\infty\}$ will
denote the maximum amount of flow that can be output, $h(w)$, if there is some
amount amount $w$ of flow into the edge: \ie, if $T \subseteq \reals^2$ is the
set of allowable flows for an edge, then
\[
    h(w) = \sup \{ t \in \reals \mid (-w, t) \in T\}.
\]
(If the set is empty, we define $h(w) =-\infty$.) In other words, $h(w)$ is
defined as the largest amount that an edge can output given an input $(-w, 0)$.
There is, of course, a natural `inverse' function which takes in an output
instead, but only one such function is necessary. Since the set $T$ is closed
by assumption, the supremum, when finite, is achieved so we have that
\[
    (-w, h(w)) \in T.
\]
We can also view $h$ as a specific parametrization of the boundary of the set
$T$ that will be useful in what follows.

\paragraph{Lossless edge example.} A simple practical example of a gain
function is the gain function for an edge which conserves flows and has finite
capacity, as in~\eqref{eq:Ti-simple-graph}:
\[
    T = \{z \in \reals^2 \mid 0 \le z_2 \le b, ~~ z_1+z_2 = 0\}.
\]
In this case, it is not hard to see that
\begin{equation}\label{eq:gain-lossless}
    h(w) = \begin{cases}
        w     & 0 \le w \le b\\
        -\infty & \text{otherwise}.
    \end{cases}
\end{equation}
The fact that $h$ is finite only when $w \ge 0$ can be interpreted as `the edge
only accepts incoming flow in one direction.'

\paragraph{Nonlinear power loss.} A more complicated example is the allowable
flow set in the optimal power flow example~\eqref{eq:app-opf-Ti}, which is,
for some convex function $\ell: \reals_+ \to \reals$,
\[
    T_i = \left\{ z \in \reals^2 \mid -b \le z_1 \le 0, ~~ z_2 \le -z_1 - \ell(-z_1) \right\}.
\]
The resulting gain function is again fairly easy to derive:
\[
    h(w) = \begin{cases}
        w - \ell(w) & 0 \le w \le b\\
        -\infty         & \text{otherwise}.
    \end{cases} 
\]
Note that, if $\ell = 0$, then we recover the original lossless edge example.
Figure~\ref{fig:opf-loss} displays this set of allowable flows $T$ and the 
associated gain function $h$.

\subsubsection{Properties}
The gain function $h$ is concave because the allowable flows set $T$ is convex,
and we can interpret the positive directional derivative of $h$ as the
current marginal price of the output flow, denominated in the input flow. Defining
this derivative as
\begin{equation}\label{eq:marginal}
    h^+(w) = \lim_{\delta \to 0^+} \frac{h(w + \delta) - h(w)}{\delta},
\end{equation}
then $h^+(0)$ is the marginal change in output if a small amount of flow were to
be added when the edge is unused, while $h^+(w)$ denotes the marginal change in
output for adding a small amount $\eps > 0$ to a flow of size $w$, for very
small $\eps$. In the case of financial markets, this derivative is sometimes
referred to as the `price impact function'. We define a reverse derivative:
\[
    h^-(w) = \lim_{\delta \to 0^+} \frac{h(w) - h(w - \delta)}{\delta},
\]
which acts in much the same way, except the limit is approached in the opposite
direction. (Both limits are well defined as they are the limits of functions
monotone on $\delta$ since $h$ is concave.) Note that, if $h$ is differentiable
at $w$, then, of course, the left and right limits are equal to the derivative,
\[
    h'(w) = h^+(w) = h^-(w),
\]
but this need not be true since we do not assume differentiability of the
function $h$. Indeed, in many standard applications, $h$ is piecewise linear and
therefore unlikely to be differentiable at optimality. On the other hand, since
the function $h$ is concave, we know that
\[
    h^+(w) \le h^-(w),
\]
for any $w \in \reals$.

\paragraph{Two-node subproblem.}
Equipped with the gain function, we can specialize the problem~\eqref{eq:subproblem-f}.
We define the arbitrage problem~\eqref{eq:subproblem-f} for an edge with gain
function $h$ as
\begin{equation}\label{eq:scalar-arb}
\begin{aligned}
& \text{maximize}   && -\eta_1w + \eta_2h(w),
\end{aligned}
\end{equation}
with variable $w \in \reals$. Since $h$ is concave, the problem is a scalar 
convex optimization problem, which can be easily solved by bisection (if the 
function $h$ is subdifferentiable) or ternary search. Since we know that
$\eta \ge 0$ by the constraints of the dual problem~\eqref{eq:dual-problem}, the
optimal value of this problem~\eqref{eq:scalar-arb} and that of the
subproblem~\eqref{eq:subproblem-f} are identical. 

\paragraph{Optimality conditions.}
The optimality conditions for problem~\eqref{eq:scalar-arb} are
that $w^\star$ is a solution if, and only if,
\begin{equation}\label{eq:swap-opt-cond}
    \eta_2 h^+(w^\star) \le \eta_1 \le \eta_2 h^-(w^\star).
\end{equation}
If the function $h$ is differentiable then $h^+ = h^- = h'$ and the expression
above simplifies to finding a root of a monotone function:
\begin{equation}\label{eq:root}
\eta_2 h'(w^\star) = \eta_1.
\end{equation}
If there is no root and condition~\eqref{eq:swap-opt-cond} does not hold, then
$w^\star = \pm\infty$. However, the solution will be finite for any feasible
flow set that does not contain a line; \ie, if the edge cannot create `infinite
flow'.

\paragraph{No-flow condition.}
The inequality~\eqref{eq:swap-opt-cond} gives us a simple way of verifying
whether we will use an edge with allowable flows $T$, given some prices $\eta_1$
and $\eta_2$. In particular, not using this edge is optimal whenever
\[
    h^+(0) \le \frac{\eta_1}{\eta_2} \le h^-(0).
\]
We can view the interval $[h^+(0), h^-(0)]$ as a `no-flow interval' for the
edge with feasible flows $T$. (In many markets, for example, this interval is a
bid-ask spread related to the fee required to place a trade.) This `no-flow
condition' lets us save potentially wasted effort of computing an optimal
arbitrage problem, as most flows in the original problem will be 0 in many
applications. In other words, an optimal flow often will not use most edges.

\subsubsection{Bounded edges} \label{sec:bounded-two-node}
In some cases, we can similarly easily check when an edge will be saturated. We
say an edge is \emph{bounded in forward flow} if there is a finite $w^0$ such
that $h(w^0) = \sup h$; \ie, there is a finite input $w^0$ which will give the
maximum possible amount of output flow. An edge is \emph{bounded} if it is
bounded in forward flow by $w^0$ and if the set $\dom h \cap (-\infty, w^0]$ is
bounded. Capacity constraints, such as those of~\eqref{eq:Ti-simple-graph}, 
imply an edge is bounded.

\paragraph{Minimum supported price.} In the dual problem, a bounded edge then
has a notion of a `minimum price'. First, define
\[
    w^\mathrm{max} = \inf \{w \in \reals \mid h(w) = \sup h\},
\]
\ie, $w^\mathrm{max}$ is the smallest amount of flow that can be tendered to maximize the
output of the provided edge. We can then define the \emph{minimum supported
price} as the left derivative of $h$ at $w^\mathrm{max}$, which is written $h^-(w^\mathrm{max})$,
from before. The first-order optimality conditions imply that $w^\mathrm{max}$ is a
solution to the scalar optimal arbitrage problem~\eqref{eq:scalar-arb} whenever
\[
    h^-(w^\mathrm{max}) \ge \frac{\eta_1}{\eta_2}.
\]
In English, this can be stated as: if the minimum supported marginal price we
receive for $w^\mathrm{max}$ is still larger than the price being arbitraged against,
$\eta_1/\eta_2$, it is optimal use all available flow in this edge.

\paragraph{Active interval.} Defining $w^\mathrm{min}$ as
\[
    w^\mathrm{min} = \inf (\dom h \cap (-\infty, w^\mathrm{max}]),
\]
where we allow $w^\mathrm{min} = -\infty$ to mean that the edge is not bounded. We then
find that the full problem~\eqref{eq:scalar-arb} needs to be solved only when
\begin{equation}\label{eq:interval}
    h^-(w^\mathrm{max}) < \frac{\eta_1}{\eta_2} < h^+(w^\mathrm{min}).
\end{equation}
We will call this interval of prices the \emph{active interval} for an edge, as
the optimization problem~\eqref{eq:scalar-arb} only needs to be solved when the
prices $\eta$ are in the interval~\eqref{eq:interval}, otherwise, the solution
is one of $w^\mathrm{min}$ or $w^\mathrm{max}$.

\subsection{Restoring primal feasibility}\label{sec:dual-primal-feas}
Unfortunately, dual decomposition methods do not, in general, find a primal
feasible solution; given optimal dual variables $\eta^\star$ and $\nu^\star$ for
the dual problem~\eqref{eq:dual-problem}, it is not the case that all solutions
$y^\star$, $x^\star$, and $\tilde x^\star$ for the
subproblems~\eqref{eq:subproblems} satisfy the constraints of the original
augmented problem~\eqref{eq:augmented-problem}. Indeed, we are guaranteed only
that \emph{some} solution to the subproblems satisfies these constraints. We
develop a second phase of the algorithm to restore primal feasibility.

For this subsection, we will assume that the net flow utility $U$ is strictly
concave, and that the edge utilities $\{V_i\}$ are each either strictly concave or
identically zero. If $V_i$ (or $U$) is nonzero, then it has a unique solution
for its corresponding subproblem at the optimal dual variables. This, in turn,
implies that the solutions to the dual subproblems are feasible, and therefore
optimal, for the primal problem. However, when some edge utilities are zero and
the corresponding sets of allowable flows are not strictly convex, we must take
care to recover edge flows that satisfy the net flow conservation constraint.

We note that, if the $\{V_i\}$ are all strictly concave (\ie, none are equal to
zero) with no restrictions on $U$, one may directly construct a solution $(y,
\{x_i\})$ by setting $x_i = \tilde x_i^\star$, the solutions to the arbitrage
subproblems for optimal dual variables $\eta^\star$ and $\nu^\star$. We can
then set
\[
    y = \sum_{i=1}^m A_i x_i,
\]
to get feasible---and therefore optimal---flows for
problem~\eqref{eq:main-problem}.

\paragraph{Example.}
Consider a lossless edge with capacity constraints, which has the allowable flow
set
\begin{equation}\label{eq:ex-bounded}
    T = \{(z_1, z_2) \mid 0 \le z_2 \le b, ~~ z_1 + z_2 = 0 \}.
\end{equation}
The associated gain function is $h(w) =  w$, if $0 \le w \le b$, and
$h(w)=-\infty$ otherwise. This gives the arbitrage
problem~\eqref{eq:scalar-arb} for the lossless edge
\[
\begin{aligned}
& \text{maximize}   && -\eta_1 w + \eta_2 w\\
& \text{subject to} && 0 \le w \le b.
\end{aligned}
\]
Proceeding analytically, we see that the optimal solutions to this problem are
\[
    w^\star \in 
    \begin{cases}
        \{0\} & \eta_1 > \eta_2 \\
        \{b\} & \eta_1 < \eta_2 \\
        [0, b]  & \eta_1 =  \eta_2.
    \end{cases}
\]
In words, we will either use the full edge capacity if $\eta_1 < \eta_2$,
or we will not use the edge if $\eta_1 > \eta_2$. However, if $\eta_1 =
\eta_2$, then any usage from zero up to capacity is an optimal solution for the
arbitrage subproblem. Unfortunately, not all of these solutions will return a
primal feasible solution for the original problem~\eqref{eq:augmented-problem}.

\paragraph{Dual optimality.}
More generally, given an optimal dual point $(\nu^\star, \eta^\star)$, an
optimal flow over edge $i$ (\ie, a flow that solves the original
problem~\eqref{eq:main-problem}) given by $x_i^\star$, will satisfy
\[
    x_i^\star \in \partial f_i(\eta_i^\star),
\]
by strong duality, as does the solution $\tilde x_i^\star$ to the arbitrage
subproblem~\eqref{eq:subproblem-f},
\[
    \tilde x_i^\star \in \partial f_i(\eta_i^\star),
\]
by definition. The subdifferential $\partial f_i(\eta_i^\star)$ is a closed
convex set, as it is the intersection of hyperplanes defined by subgradients.
We distinguish between two cases. First, if the set $T_i$ is strictly convex,
then the set $\partial f_i(\eta_i^\star)$ consists of a single point and
$x_i^\star = \tilde x_i^\star$. However, if $T_i$ is not strictly convex, then
we only are guaranteed that
\[
        x_i^\star \in T_i^\star(\eta^\star_i) = 
        T_i \cap \partial f_i(\eta_i^\star).
\]
This set $T^\star_i(\eta_i^\star)$ is the intersection of two convex sets and,
therefore, is convex. In fact, this set is exactly a `face' of $T_i$ with
supporting hyperplane defined by $\eta_i^\star$. In general, this set can be as
hard to describe as the original set $T_i$. On the other hand, in the common
case that $T_i$ is polyhedral or two-dimensional, the set has a
concise representation that is easier to optimize over than the set itself.
(Note that, in practice, numerical precision issues may also need to be taken
into account, as we only know $\eta_i^\star$ up to some tolerance.)


\paragraph{Two-node edges.}
For two-node edges, observe that a piecewise linear set of allowable flows
$T_i$ can be written as a Minkowski sum of its segments. Equivalently, a
piecewise linear gain function is equivalent to adding bounded linear edges for
each of its segments (\cf~\eqref{eq:Ti-simple-graph}). For a given optimal
price vector $\eta_i^\star$, the optimal flow $x_i^\star$ will be nonzero on at
most one of these segments, and the set $T_i^\star$ is a single point unless
$\eta_i^\star$ is normal to one of these line segments. This idea, of course,
may be extended to general two-node allowable flows whose boundary may include
smooth regions as well as line segments. Returning to
example~\eqref{eq:ex-bounded} above, if $\eta_i^\star = \alpha\ones$ for some
$\alpha > 0$, then
\[
    T_i^\star(\eta_i^\star) = \{z \in \reals^2 \mid \ones^T z = 0, ~~ 0 \le z_2 \le b\}.
\]
Otherwise, $T_i^\star(\eta_i^\star)$ is an endpoint of this line segment:
either
\[
    T_i^\star(\eta_i^\star) = \{(0, 0)\},
\]
or
\[
    T_i^\star(\eta_i^\star) = \{(-b, b)\}.
\]

\paragraph{Recovering the primal variables.}
Recall that the objective function $U$ is strictly concave by assumption, so
there is a unique solution that solves the associated
subproblem~\eqref{eq:subproblem-U} at optimality. Let $S$ be a set containing
the indices of the strictly convex feasible flows; that is, the index $i \in S$
if $T_i$ is strictly convex. Now, let the dual optimal points be $(\nu^\star,
\eta^\star)$, and the optimal points for the
subproblems~\eqref{eq:subproblem-U} and~\eqref{eq:subproblem-f} be $y^\star$
and $\tilde x_i^\star$ respectively. We can then recover the primal variables
by solving the problem
\[
\begin{aligned}
    &\text{minimize} && {\textstyle \|y^\star - \sum_{i=1}^m A_i x_i\| }\\
    &\text{subject to} && x_i = \tilde x_i^\star, \quad i \in S \\
    &&& x_i \in T_i^\star(\eta^\star_i), \quad i \notin S.
\end{aligned}
\]
Here, the objective is to simply find a set of feasible $x_i$ (\ie, that `add
up' to $y^\star$) which are consistent with the dual prices discovered by the
original problem, in the sense that they minimize the error between their net
flows and the net flow vector $y^\star$. Indeed, if the problem is correctly
specified (and solution errors are not too large), the optimal value should
always be 0. When the sets $\{T_i^\star\}$ can be described by linear
constraints and we use the $\ell_1$ or $\ell_\infty$ norm, this problem is a
linear program and can be solved very efficiently. The two-node linear case is
most common in practice, and we leave further exploration of the reconstruction
problem to future work.

\section{Numerical examples}\label{sec:numerical-examples}
We illustrate our interface by revisiting some of the examples
in~\S\ref{sec:apps}. We do not focus on the linear case, as this case is better
solved with special-purpose algorithms such as the network simplex method or the
augmenting path algorithm. In all experiments, we examine the convergence of our
method, \texttt{ConvexFlows}, and compare its runtime to the commercial solver
Mosek~\cite{mosek}, accessed through the JuMP modeling 
language~\cite{dunningJuMPModelingLanguage2017,legatMathOptInterfaceDataStructure2021,Lubin2023}.
We note that the conic formulations of these problems often do not preserve the
network structure and may introduce a large number of additional variables.

Our method \texttt{ConvexFlows} is implemented in the Julia programming language~\cite{bezansonJuliaFreshApproach2017},
and the package may be downloaded from
\begin{center}
    \texttt{\url{https://github.com/tjdiamandis/ConvexFlows.jl}}.
\end{center}
Code for all experiments is in the \texttt{paper} directory of the repository.
All experiments were run using \texttt{ConvexFlows} v0.1.1 on a  MacBook Pro with 
a M1 Max processor (8 performance cores) and 64GB of RAM. We suspect that our
method could take advantage of further parallelization than what is available
on this machine, but we leave this for future work.

\subsection{Optimal power flow}
We first consider the optimal power flow problem from~\S\ref{sec:app-opf}. This
problem has edges with only two adjacent nodes, but each edge flow has a
strictly concave gain function due to transmission line loss. These line losses
are given by the constraint set~\eqref{eq:app-opf-Ti}, and we use the objective
function~\eqref{eq:app-opf-U} with the quadratic power generation cost functions
\[
    c_i(w) = \begin{cases}
        (1/2)w^2 & w \ge 0 \\
        0 & w < 0.
    \end{cases}
\]
Since the flow cost functions are identically zero, we only have two subproblems
(\cf,~\S\ref{sec:special-cases}). The first subproblem is the evaluation of
$\bar U$, which can be worked out in closed form:
\[
    \bar U(\nu) = (1/2)\|\nu\|_2^2 - d^T \nu,
\]
with domain $\nu \ge 0$. We could easily add additional constraints, such as an
upper bound on power generation, but do not for the sake of simplicity. The
second subproblem is the arbitrage problem~\eqref{eq:scalar-arb},
\[
    f_i(\eta_i) = \sup_{0 \le w \le b_i} \left\{
         -\eta_1 w  + \eta_2 \left( w - \ell_i(w) \right)
    \right\},
\]
which can generally be solved as a single-variable root finding problem because
the allowable flows set is strictly convex. Here, the edge is, in addition, 
`bounded' (see~\S\ref{sec:dual-two-node}) with a closed form solution. We
provide the details in appendix~\ref{app:opf}.

\begin{figure}[h]
    \centering
    \includegraphics[width=0.5\textwidth]{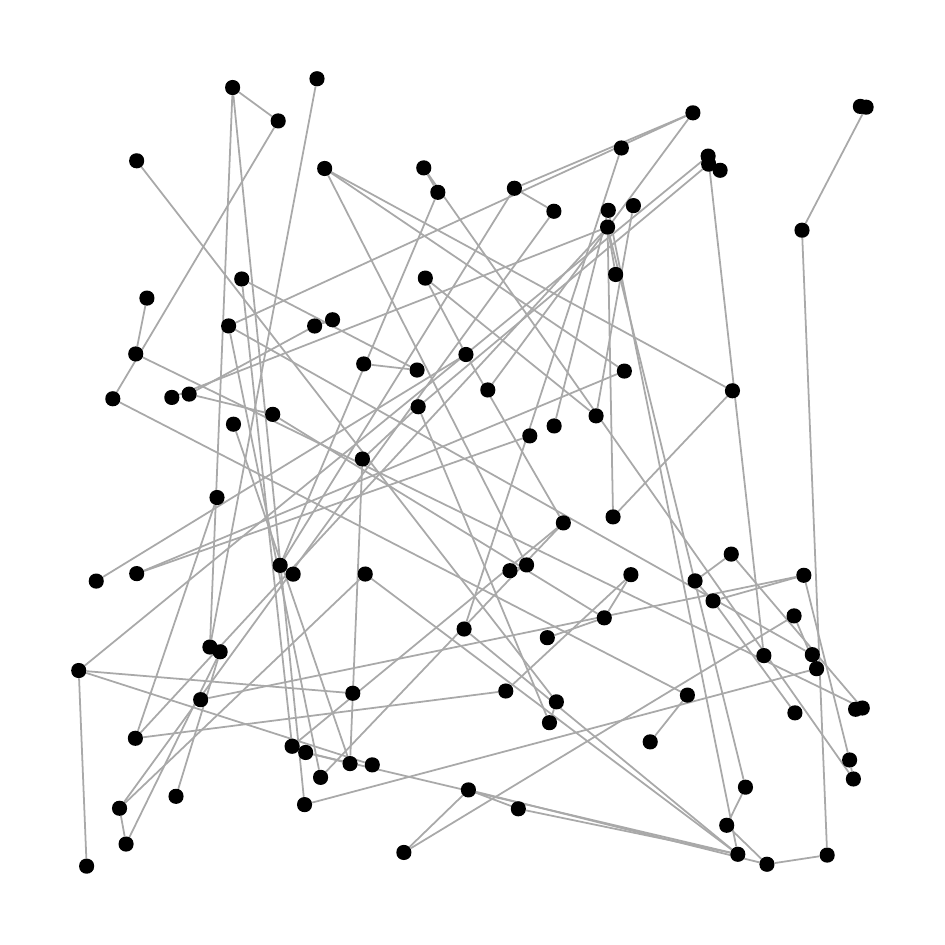}
    \caption{Sample network for $n = 100$ nodes.}
    \label{fig:opf-network}
\end{figure}

\paragraph{Problem data.} 
We model the network as in~\cite{kraning2013dynamic} with the same parameters,
which results in a network with high local connectivity and a few longer
transmission lines. Figure~\ref{fig:opf-network} shows an example with $n =
100$. We draw the demand $d_i$ for each node uniformly at random from the set
$\{0.5, 1, 2\}$. For each transmission line, we set $\alpha_i = 16$ and $\beta_i
= 1/4$. We draw the maximum capacity for each line uniformly at random from the
set $\{1, 2, 3\}$. These numbers imply that a line with maximum capacity $1$
operating at full capacity will loose about $10\%$ of the power transmitted,
whereas a line with maximum capacity $3$ will loose almost $40\%$ of the power
transmitted. For the purposes of this example, we let all lines be
bidirectional: if there is a line connecting node $j$ to node $j'$, we add a
line connecting node $j'$ to node $j$ with the same parameters.

\begin{figure}[h]
    \begin{subfigure}[t]{0.46\textwidth}
        \centering
        \includegraphics[width=\columnwidth]{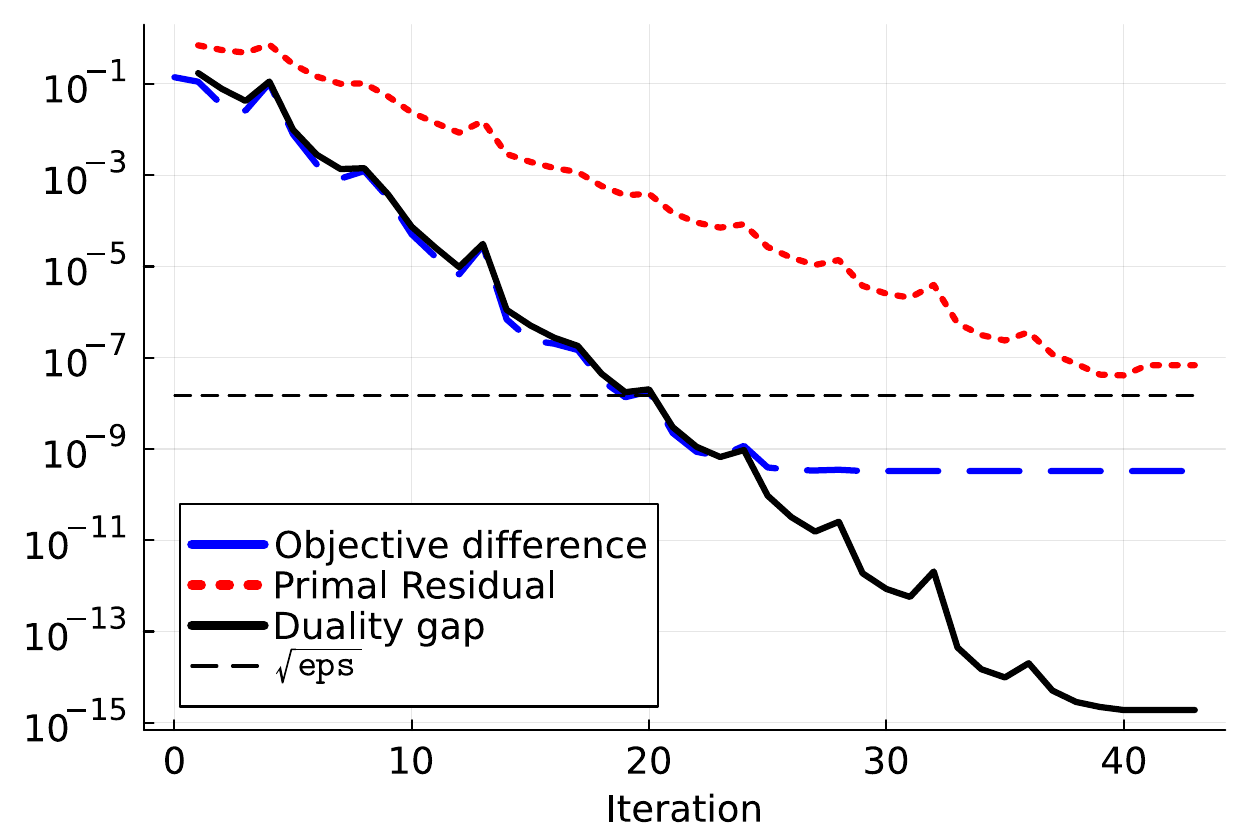}
        \caption{
            Convergence of \texttt{ConvexFlows} with $n = 100$. The objective is
            compared to a high-precision solution from Mosek. The primal
            residual measures the net flow constraint violation, with $\{x_i\}$
            from~\eqref{eq:subproblem-f} and $y$ from~\eqref{eq:subproblem-U}.
        }
        \label{fig:opf-numerics-iter}
    \end{subfigure}
    \hfill
    \begin{subfigure}[t]{0.46\textwidth}
        \centering
        \includegraphics[width=\columnwidth]{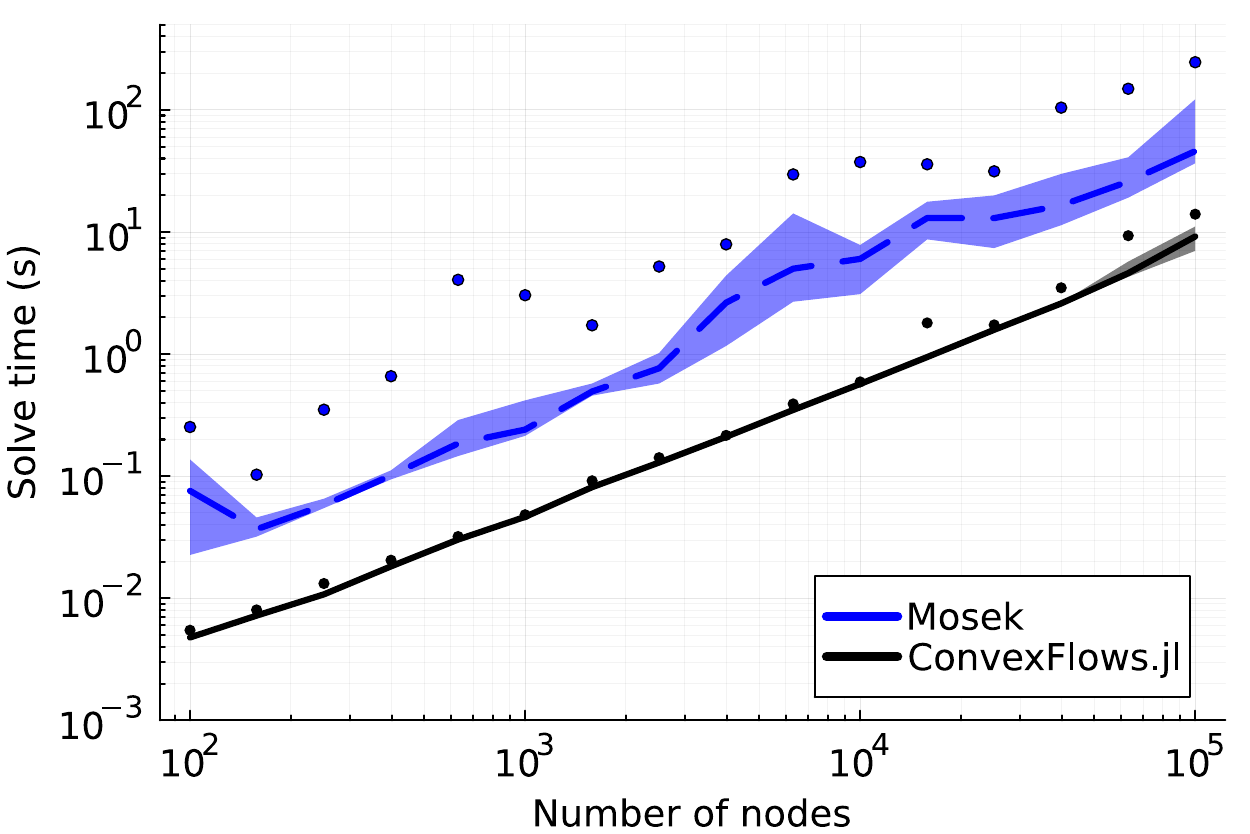}
        \caption{
            Comparison of \texttt{ConvexFlows} with Mosek. Lines indicate the
            median time over 10 trials, and the shaded region indicates the 25th
            to 75th quantile range. Dots indicate the maximum time over the 10
            trials.
        }
        \label{fig:opf-numerics-compare}
    \end{subfigure}
    \caption{Numerical results for the optimal power flow problem.}
    \label{fig:opf-numerics}
\end{figure}

\paragraph{Numerical results.}
We first examine the convergence per of our method on an example with $n = 100$
nodes and $m = 198$ transmission lines. In figure~\ref{fig:opf-numerics-iter},
we plot the relative duality gap, net flow constraint violation, and difference 
between our objective value and the `optimal' objective value, obtained using
the commercial solver Mosek. (See appendix~\ref{app:opf} for the conic
formulation). These results suggest that our method enjoys linear
convergence. The difference in objective value at `optimality' is likely due to
floating point numerical inaccuracies, as it is below the square root of machine
precision, denoted by $\sqrt{\texttt{eps}}$. In figure~\ref{fig:opf-numerics-compare}, we
compare the runtime of our method to Mosek for increasing values of $n$, with
ten trials for each value. For each $n$, we plot the median time, the 25th to
75th quantile, and the maximum time. Our method clearly results in a significant
and consistent speedup. Notably, our method exhibits less variance in solution
time as well. We emphasize, however, that our implementation is not highly
optimized and relies on an `off-the-shelf' L-BFGS-B solver. We expect that
further software improvement could yield even better performance.

\subsection{Routing orders through financial exchanges}
Next, we consider a problem which includes both edges connecting more than 
two nodes and utilities on the edge flows: routing trades through decentralized 
exchanges (see~\S\ref{sec:app-cfmm}). For all experiments, we use
the net flow utility function 
\[
    U(y) = c^Ty - I_{\reals_+^n}(y).
\]
We interpret this function as finding arbitrage in the network of markets. More
specifically, we wish to find the most valuable trade $y$, measured according to
price vector $c$, which, on net, tenders no assets to the network. The associated
subproblem~\eqref{eq:subproblem-U} can be computed as
\[
    \bar U(\nu) = \begin{cases}
        0 & \nu \ge c \\
        \infty & \text{otherwise}.
    \end{cases}
\]
We also want to ensure our trade with any one market is not too large, so we
add a penalty term to the objective:
\[
    V_i(x_i) = -(1/2)\|(x_i)_-\|_2^2.
\]
(Recall that negative entries denote assets tendered to an exchange.) The
associated subproblem is
\[
    \bar V_i(\xi) = (1/2) \|\xi\|_2^2.
\]
Finally, the arbitrage problem here is exactly the problem of computing an
optimal arbitrage trade with each market, given prices on some infinitely-liquid
external market.

\paragraph{Problem data.}
We generate $m$ markets and $n = 2 \sqrt{m}$ assets. We use three popular market
implementations for testing: Uniswap (v2)~\cite{uniswap},
Balancer~\cite{balancer} two-asset markets, and Balancer three-asset markets.
We provide the details of these markets and the associated arbitrage problems in
appendix~\ref{app:cfmm}. Markets are Uniswap-like with probability $2/5$,
Balancer-like two-asset markets with probability $2/5$, and Balancer-like
three-asset markets with probability $1/5$. Each market $i$ connects randomly
selected assets and has reserves sampled uniformly at random from the interval
$[100,\, 200]^{n_i}$.

\begin{figure}
    \begin{subfigure}[t]{0.46\textwidth}
        \centering
        \includegraphics[width=\columnwidth]{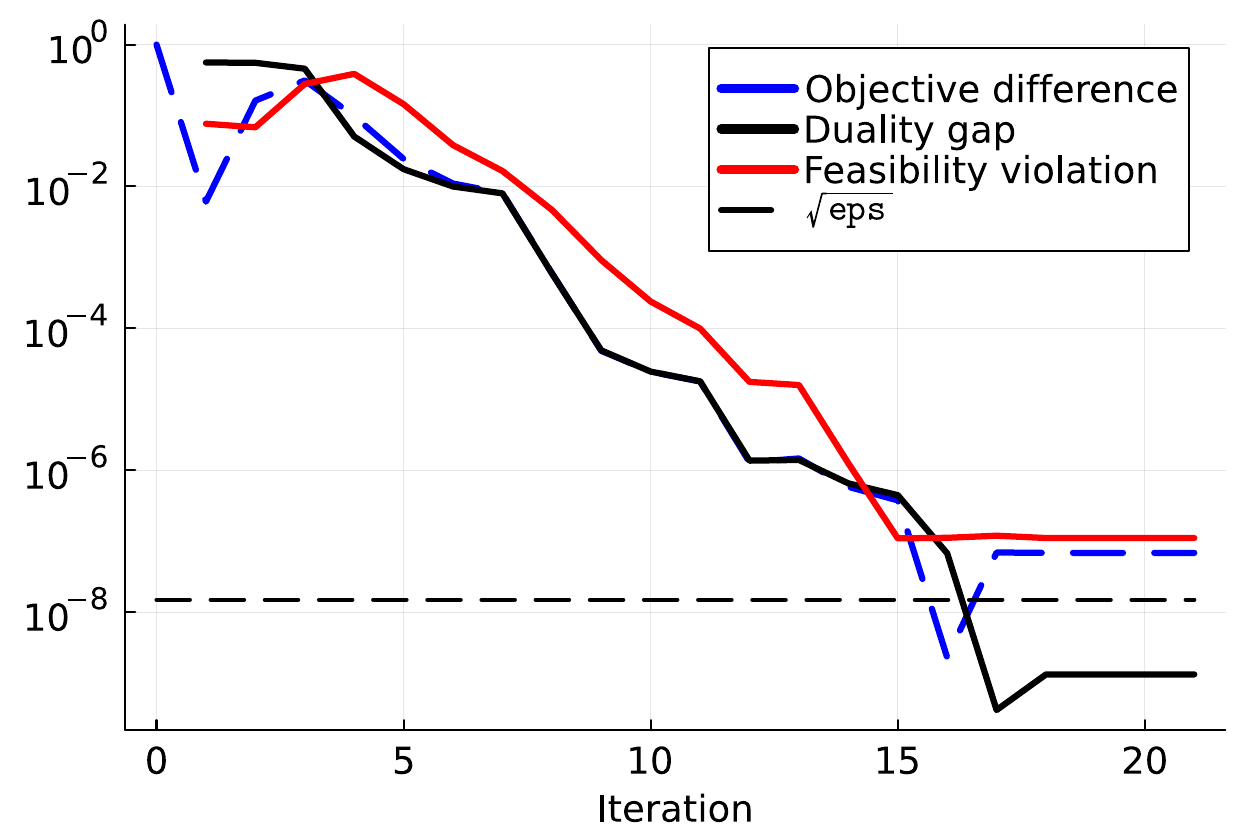}
        \caption{Convergence without edge penalties.}
        \label{fig:cfmm-numerics-iter}
    \end{subfigure}
    \hfill
    \begin{subfigure}[t]{0.46\textwidth}
        \centering
        \includegraphics[width=\columnwidth]{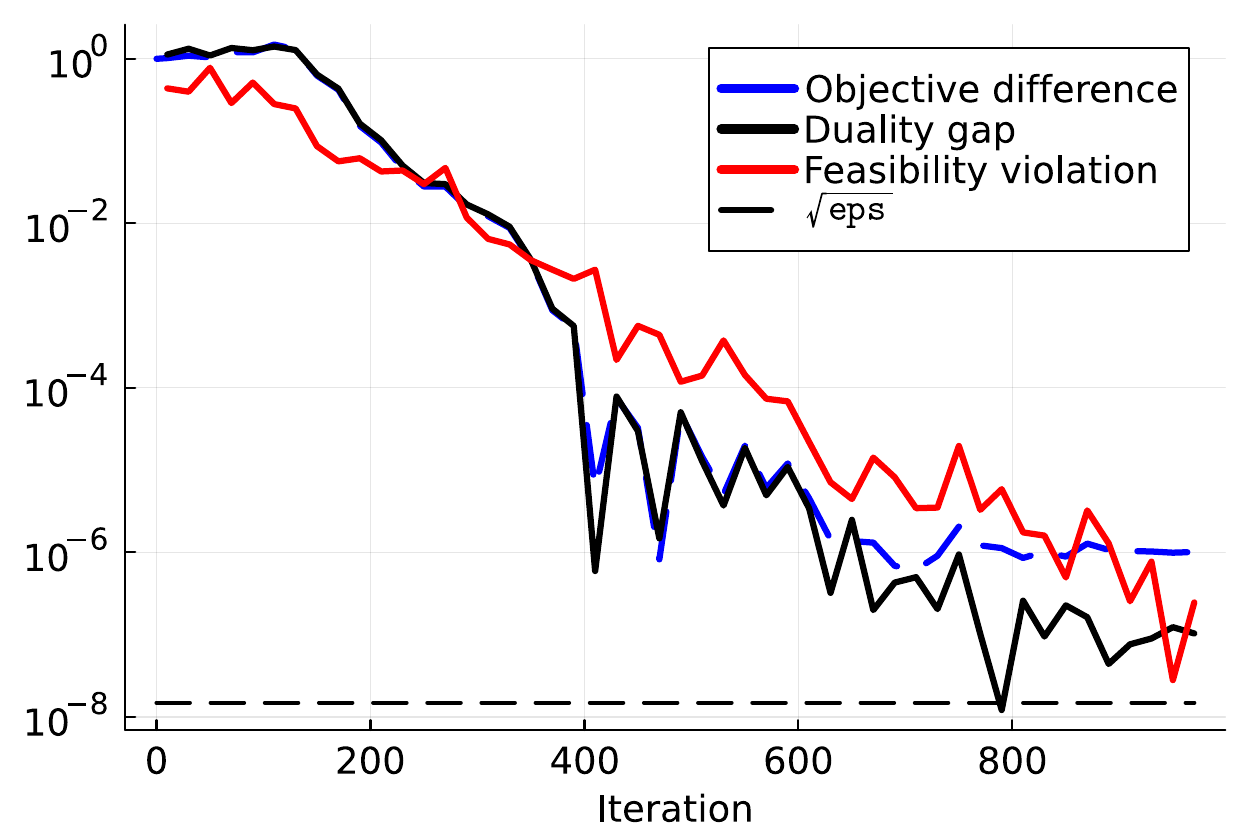}
        \caption{Convergence with edge penalties.}
        \label{fig:cfmm-numerics-iter-vi}
    \end{subfigure}
    \caption{Convergence of \texttt{ConvexFlows} on an example with $n = 100$ 
    assets and $m = 2500$ markets.}
    \label{fig:cfmm-numerics}
\end{figure}

\begin{figure}
    \centering
    \includegraphics[width=0.6\textwidth]{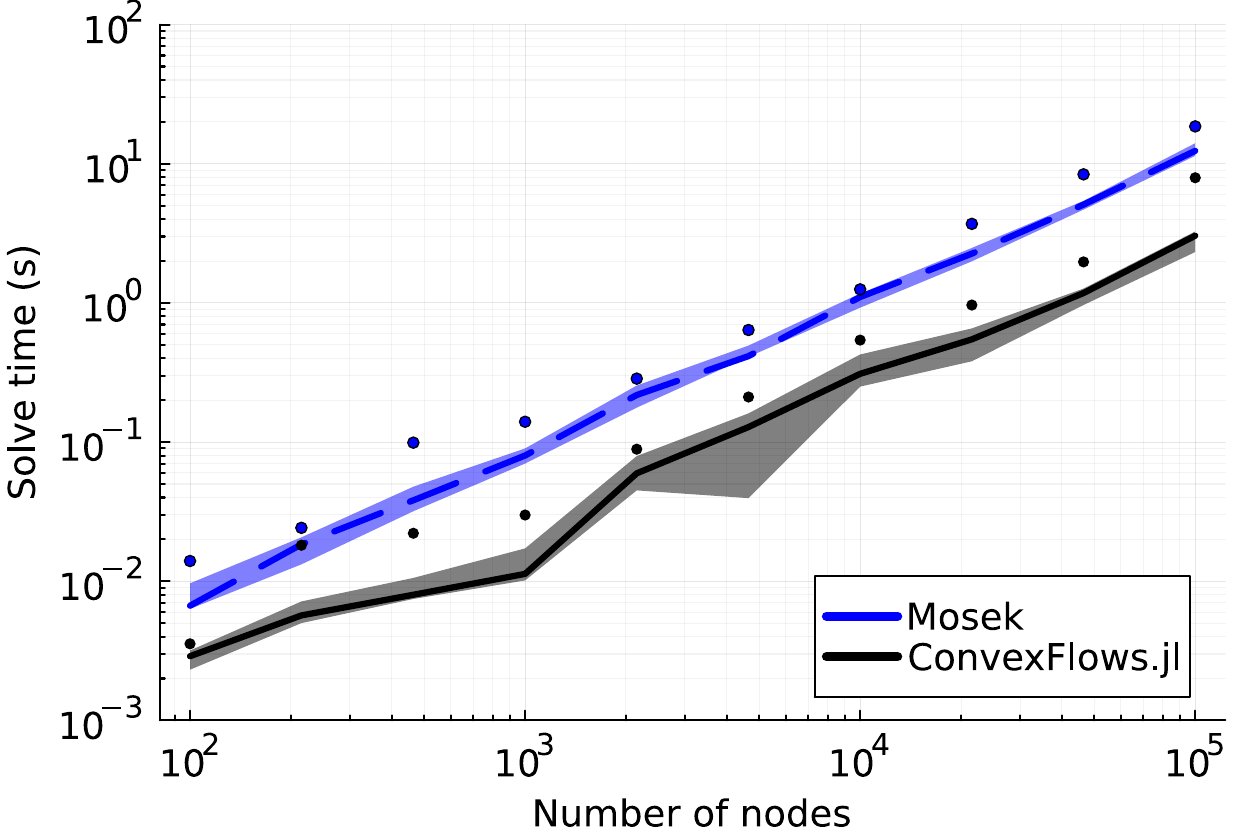}
    \caption{Comparison of \texttt{ConvexFlows} and Mosek for $m$ varying from
    $100$ to $100,000$ and $n = 2\sqrt{m}$. Lines indicate the median time over
    10 trials, and the shaded region indicates the 25th to 75th quantile range.
    Dots indicate the maximum time over the 10 trials.}
    \label{fig:cfmm-numerics-compare}
\end{figure}

\paragraph{Numerical results.}
We first examine the convergence of our method on an example with $m = 2500$ and
$n = 100$. In figures~\ref{fig:cfmm-numerics-iter}
and~\ref{fig:cfmm-numerics-iter-vi}, we plot the convergence of the relative
duality gap, the feasibility violation of $y$, and the relative difference
between the current objective value and the optimal objective value, obtained
using the commercial solver Mosek. (See appendix~\ref{app:cfmm} for the conic
formulation we used.) Note that, here, we reconstruct $y$ as
\[
    y = \sum_{i=1}^m A_i x_i,
\]
instead of using the solution to the subproblem~\eqref{eq:subproblem-U} as we
did in the previous example. As a result, this $y$ satisfies the net flow
constraint by construction. We measure the feasibility violation relative to the
implicit constraint in the objective function $U$, which is that $y \ge 0$. We
again see that our method enjoys linear convergence in both cases; however, the
convergence is significantly slower when edge objectives are added
(figure~\ref{fig:cfmm-numerics-iter-vi}). We then compare the runtime of our
method and the commercial solver Mosek, both without edge penalties and with
only two-node edges, in figure~\ref{fig:cfmm-numerics-compare}. Again,
\texttt{ConvexFlows} enjoys a significant speedup over Mosek.

\section{Conclusion}
In this paper, we introduced the convex network flow problem, which is a 
natural generalization of many important problems in computer science, 
operations research, and related fields. We showed that many problems from the
literature are special cases of this framework, including max-flow, optimal 
power flow, routing through financial markets, and equilibrium computation in
Fisher markets, among others. This generalization has a number of useful 
properties including, and perhaps most importantly, that its dual decomposes 
over the (hyper)graph structure. This decomposition results in a fast algorithm 
that easily parallelizes over the edges of the graph and preserves the
structure present in the original problem. We implemented this algorithm in the
Julia package \texttt{ConvexFlows.jl} and showed order-of-magnitude speedups over
a commercial solver applied to the same problem.

\paragraph{Future work.}
We believe that analyzing the convex flow problem properties in a more systematic
way will lead to several interesting future research directions. For example, 
we mention in~\S\ref{sec:app-max-flow} that bidirectional edge flows can be 
viewed as the Minkowski sum of two directional edge flows, one in each direction. 
Can this idea be generalized to other feasible sets? What does a natural version
of this look like? Other important generalizations include: is it possible to
extend this framework to include fixed costs for using (or not using) an edge?
Does this problem's dual formulation have a more natural dual that has a
similar interpretation as the primal, akin to the self-duality in extended
monotropic programming~\cite{bertsekas2008extended}? And, finally, is there an
easier-to-use interface for solvers of this particular problem, which does not
require specifying solutions to the subproblems~\eqref{eq:subproblems}
directly? We suspect many of these questions are interesting from both a
theoretical and practical perspective, given the relevance of this problem
formulation to many applications.

\ifsubmit
    \begin{acknowledgements}
\else
    \section*{Acknowledgements}
\fi
The authors thank Flemming Holtorf, Pablo Parrilo, and Anthony Degleris for
helpful discussions. Theo Diamandis is supported by the Department of Defense
(DoD) through the National Defense Science \& Engineering Graduate (NDSEG)
Fellowship Program. This material is based upon work supported by the National
Science Foundation under Award No.\ OSI-2029670. Any opinions, findings and
conclusions or recommendations expressed in this material are those of the
authors and do not necessarily reflect the views of the National Science
Foundation.
\ifsubmit
    \end{acknowledgements}
\else\fi

\ifsubmit
    \bibliography{refs}
\else
    \printbibliography
\fi

\appendix

\section{Extended monotropic programming}~\label{app:extended-monotropic}
In this section, we explicitly draw the connection between the extended
monotropic programming (EMP) problem formulated by Bertsekas~\cite{bertsekas2008extended}
and the convex flow problem~\eqref{eq:main-problem}. The extended 
monotropic programming problem can be written as
\[
\begin{aligned}
    &\text{minimize} && \sum_{i=1}^{m+1} f_i(z_i) \\
    &\text{subject to} && z \in S,
\end{aligned}
\]
with variable $z \in \reals^N$. The functions $f_i$ are convex functions of the subvectors 
$z_i$, and the set $S$ is a subspace of $\reals^N$. Taking $z = (y, x_1, \dots, x_m)$
and changing the minimization to a maximization, we can write the convex 
flow problem as a monotropic programming problem:
\[
    \begin{aligned}
        &\text{maximize} && U(y) + \sum_{i=1}^{m} V_i(x_i) - I_{T_i}(x_i)  \\
        &\text{subject to} && y = \sum_{i=1}^m A_i x_i,
    \end{aligned}
\]
where we took
\[
    f_{m+1} = -U, \qquad \text{and} \qquad f_i = -V_i + I_{T_i}, \quad i = 1, \dots, m,
\]
Note that the linear net flow constraint is a subspace constraint.

\paragraph{Duality.}
The dual of the EMP problem considered by Bertsekas
is given by
\[
\begin{aligned}
    &\text{maximize} && -\sum_{i=1}^{m+1} \sup_{z_i \in \reals^{n_i}} \left\{ \lambda_i^Tz_i - f_i(z_i) \right\} \\
    &\text{subject to} && \lambda \in S^\perp.
\end{aligned}
\]
Substituting in $U$ and $\{V_i\}$ and switching the sign of $\lambda$, the 
objective terms become
\[
    \sup_{z_{m+1} \in \reals^n}\left\{ U(z_{m+1}) - \lambda_{m+1}^Tz_{m+1} \right\} = \bar U(\lambda_{m+1})
    \qquad \text{and} \qquad
    \sup_{z_i \in T_i}\left\{ V_i(z_i) - \lambda_i^Tz_i \right\}.
\]
These terms are very close to, but not exactly the same as, the dual terms in
the convex flow problem~\eqref{eq:dual-problem-impl}. In particular, the
$U$ subproblem~\eqref{eq:subproblem-U} remains the same, but, in our framework,
we introduced an additional dual variable to split the $V_i$ subproblem into two
subproblems: one for the function $V_i$~\eqref{eq:subproblem-V} and one for the
set $T_i$~\eqref{eq:subproblem-f}. This split allows for a more efficient
algorithm that uses the `arbitrage' primitive~\eqref{eq:subproblem-f}, which has
a very fast implementation for many edges, especially in the case of two node 
edges (see~\S\ref{sec:dual-two-node}). Our dual problem for the convex flow
problem allows us to exploit more structure in our solver.

\paragraph{When the EMP problem matches.}
In the case of zero edge utilities, however, the EMP problem matches the convex
network flow problem exactly. In this case, the $V$ subproblem disappears and we
are left only with the arbitrage subproblem:
\[
    \sup_{z_i \in T_i}\left\{ V_i(z_i) - \lambda_i^Tz_i \right\} =
    \sup_{z_i \in T_i}\left\{ -\lambda_i^Tz_i \right\} = f_i(-\lambda_i).
\]
Letting $\nu = -\lambda_m$, the subspace constraint then becomes
\[
    \lambda_i = A_i^T\nu, \quad i = 1, \dots, m.
\]
Thus, we recover the exact dual of the convex flow problem with zero
edge utilities, given in~\eqref{eq:zero-edge-problem}. This immediately implies
the strong duality result given in~\cite[Prop 2.1]{bertsekas2008extended} holds
in our setting as well.

\paragraph{Self duality.}
The EMP dual problem has the same form as the primal; in this sense, the EMP
problem is self-dual. The convex flow problem, however, does not appear
to be self-dual in the same sense, since we consider a very specific subspace
that defines the net flow constraint. We leave exploration of duality in our
setting to future work.

\section{Fisher market problem KKT conditions}\label{app:fisher-market}
The Lagrangian of the Fisher market problem~\eqref{eq:app-market-equilibrium} is
\[
    L(x, \mu, \lambda) = 
    \sum_{i =1}^{n_b} b_i \log(U(x_i)) + 
    \mu^T\left(\ones - \sum_{i = 1}^{n_b} x_i\right) - 
    \sum_{i = 1}^{n_b} x_i^T\lambda_i,
\]
where $\{x_i \in \reals^{n_g}\}$ are the primal variables and 
$\mu \in \reals^{n_g}$ and $\{\lambda_i \in \reals_+^{n_g}\}$ are the dual variables.
Let $x^\star$, $\mu^\star$, $\lambda^\star$ be a primal-dual solution to this
problem.
The optimality conditions~\cite[\S5.5]{cvxbook} are primal feasibility, 
complementary slackness, and the dual condition
\[
    \partial_{x_i}L(x^\star, \mu^\star, \lambda^\star) = 
    \frac{b_i}{U(x_i^\star)}\nabla U(x_i^\star) - \mu^\star - \lambda^\star_i = 0, 
    \quad \text{for} ~ i = 1, \dots, n_b.
\]
This condition simplifies to
\[
    \nabla U(x_i^\star) \ge (U(x_i) / b_i) \cdot \mu^\star, \quad \text{for} ~ i = 1, \dots, n_b.
\]
If we let the prices of the goods be $\mu^\star \in \reals^{n_g}$, this condition says that the
marginal utility gained by an agent $i$ from an additional small amount of any 
good is at least as large as that agent's budget-weighted price times their
current utility. As a result, the prices $\mu^\star$ will cause all agents to 
spend their entire budget on a utility-maximizing basket of goods, and all goods 
will be sold.

\section{Additional details for the numerical experiments}
\subsection{Optimal power flow}\label{app:opf}

\paragraph{Arbitrage problem.}
Here, we explicitly work out the arbitrage subproblem for the optimal power flow
problem. Recall that the set of allowable flows is given by (dropping the edge
index for convenience)
\[
    T = \{z \in \reals^2 \mid -b \le z_1 \le 0 ~ z_2 \le -z_1 - \ell(-z_1)\},
\]
where 
\[
    \ell(w) = 16 \left(\log(1 + \exp(w/4)) - \log 2\right) - 2w
\]
Given an edge input $w \in [0, b]$, the gain function is
\[
    h(w) = w - \ell(w),
\]
where we assume the edge capacity $b$ is chosen such that the
function $f$ is increasing for all $w \in [0, b]$, \ie, $f'(b) > 0$.
Using~\eqref{eq:root}, we can compute the optimal solution $x^\star$ to the
arbitrage subproblem~\eqref{eq:subproblem-f} as
\[
    x^\star_1 = \left(4 \log\left(\frac{3\eta_2 - \eta_1}{\eta_2 + \eta_1}\right)\right)_{[0, b]}, \qquad
    x^\star_2 = h(x^\star_1),
\]
where $(\cdot)_{[0, b]}$ denotes the projection onto the interval $[0, b]$.

\paragraph{Conic formulation.}
Define the exponential cone as
\[
    K_\mathrm{exp} =
    \{
        (x,y,z) \in \reals^3 \mid y > 0, \; y e^{x/y} \le z.
    \}
\]
The transmission line constraint is of the form
\[
    \log(1 + e^s) \le t,
\]
which can be written as~\cite[\S5.2.5]{mosekcookbook}
\[
\begin{aligned}
    u + v &\le 1 \\
    (x-t,\, 1,\, u) &\in K_\mathrm{exp} \\
    (-t,\, 1,\, v) &\in K_\mathrm{exp}.
\end{aligned}
\]
Define the rotated second order cone as
\[
    K_\mathrm{rot2} = \{ (t,u,x) \in \reals_+ \times \reals_+ \times \reals^{n} \mid
     2tu \ge \| x \|_2^2\}.
\]
We can write the cost function
\[
    c_i(w) = (1/2)w_+^2,
\]
where $w_+ = \max(w, 0)$ denotes the negative part of $w$, in conic form as
minimizing $t_1 \in \reals$ subject to the second-order cone constraint~\cite[\S3.2.2]{mosekcookbook}
\[
    (0.5,\, t_1,\, t_2) \in K_\mathrm{rot2}, \qquad t_2 \ge w, \qquad t_2 \ge 0.
\]
Putting this together, the conic form problem is
\[
\begin{aligned}
&\text{maximize} && -\ones^T t_1 \\
&\text{subject to} && 
(0.5,\, (t_1)_i,\, (t_2)_i) \in K_\mathrm{rot2}, \quad \text{for}~ i = 1, \dots n\\
&&&t_2 \ge d - y, \qquad t_2 \ge 0 \\
&&&-b_i \le (x_i)_1 \le 0, \quad \text{for}~ i = 1, \dots m \\
&&&u_i + v_i \le 1 \quad\text{for}~ i = 1, \dots m\\
&&&\left(-\beta_i (x_i)_1 + (3(x_i)_1 + (x_i)_2)/\alpha - \log(2),\, 1,\, u_i\right) \in K_\mathrm{exp} \quad\text{for}~ i = 1, \dots m\\
&&&\left((3(x_i)_1 + (x_i)_2)/\alpha - \log(2),\, 1,\, v_i\right) \in K_\mathrm{exp} \quad\text{for}~ i = 1, \dots m.
\end{aligned}
\]

\subsection{Routing orders through financial exchanges}\label{app:cfmm}
In this example, we considered three different types of decentralized exchange
markets: Uniswap-like, Balancer-like swap markets, and Balancer-like multi-asset 
markets. Recall that a constant
function market maker (CFMM) allows trades between the $n$ tokens in its reserves
$R \in \reals^n_+$ with behavior governed by a trading function $\phi:\reals^n_+ 
\to \reals$. The CFMM only accepts a trade $(\Delta, \Lambda)$ where $\Delta \in 
\reals^n_+$ is the basket of tendered tokens and $\Lambda \in \reals^n_+$ is the 
basket of received tokens if 
\[
    \phi(R + \gamma\Delta - \Lambda) \ge \phi(R).
\]
The Uniswap trading
function $\phi_\mathrm{Uni}:\reals^2_+ \to \reals$ is given by
\[
    \phi_\mathrm{Uni}(R) = \sqrt{R_1R_2}.
\]
The Balancer swap market trading function $\phi_\mathrm{Bal}:\reals^2_+ \to 
\reals$ is given by
\[
    \phi_\mathrm{Bal}(R) = R_1^{4/5}R_2^{1/5}.
\]
The Balancer multi-asset trading function $\phi_\mathrm{Mul}:\reals^3_+ \to 
\reals$ is given by
\[
    \phi_\mathrm{Mul}(R) = R_1^{1/3}R_2^{1/3}R_3^{1/3}.
\] 
These functions are easily recognized as (weighted) geometric means and can be
verified as concave, nondecreasing function. Thus, the set of allowable trades
\[
    T = \{\Lambda - \Delta \mid \Lambda,\Delta \in \reals^n_+ ~\text{and}~
    \phi(R + \gamma\Delta - \Lambda) \ge \phi(R)\},
\]
is convex. Furthermore, the arbitrage problem~\eqref{eq:subproblem-f} has a
closed form solution for the case of the swap markets (see~\cite[App. A]{angerisAnalysisUniswapMarkets2020}
and the implementation from~\cite{diamandis2023efficient}). Multiasset pools may
have closed form solutions as well, which we discuss in the next section.

\paragraph{Separable CFMM arbitrage problem.}
Consider a separable CFMM with trading function $\phi:\reals^n_+ \to \reals$ of
the form
\[
    \phi(R) = \sum_{i=1}^n \phi_i(R_i),
\]
where each $\phi_i$ is strictly concave and increasing. (The non-strict case
follows from the same argument but requires more care.) Note that many CFMMs may
be transformed into this form. For example, a weighted geometric mean CFMM like
Balancer can be written in this form using a log transform:
\[
    \prod_{i=1}^n R_i^{w_i} \ge k \iff \sum_{i=1}^n w_i \log R_i \ge \log k.
\]
The arbitrage subproblem can be written as
\[
\begin{aligned}
    &\text{maximize} && \eta^T(\Lambda - \Delta) \\
    &\text{subject to} && \sum_{i=1}^{n} \phi_i(R_i + \gamma \Delta_i - \Lambda_i) \ge k \\
    &&& \Delta, \Lambda \ge 0.
\end{aligned}
\]
After pulling the nonnegativity constraints into the objective, the Lagrangian
is separable and can be written as
\[
    L(\Delta, \Lambda, \lambda) = 
    \sum_{i=1}^{n} \eta_i (\Lambda_i - \Delta_i) - I(\Delta_i) - I(\Lambda_i) + \lambda(\phi_i(R_i + \gamma \Delta_i - \Lambda_i) - k),
\]
where $\lambda \ge 0$ and $I$ is the nonnegative indicator function satisfying
$I(w) = 0$ if $w \ge 0$ and $+\infty$ otherwise. Maximizing over the primal
variables $\Delta$ and $\Lambda$ gives the dual function:
\begin{equation}\label{eq:app-dual-sep}
    g(\lambda) = \sum_{i=1}^n \left(\sup_{\Delta_i,\,\Lambda_i \ge 0} 
        \eta_i(\Lambda_i - \Delta_i) + \lambda_i\phi_i(R_i + \gamma \Delta_i - \Lambda_i)
    \right) - \lambda k.
\end{equation}
Consider subproblem $i$ inside of the sum. If $0 \le \gamma < 1$, then at most
one of $\Delta_i^\star$ or $\Lambda_i^\star$ is nonzero, which turns this two
variable problem into two single variable convex optimization problems, each
with a nonnegativity constraint. (This follows
from~\cite[\S2.2]{angerisConstantFunctionMarket2021}.) In particular, to solve
the original, we can solve two (smaller) problems by considering the two possible
cases. In the first case, we have $\Lambda_i = 0$ and $\Delta_i \ge 0$, giving
the problem
\[
    \begin{aligned}
        & \text{maximize} && -\eta_i \Delta_i + \lambda_i \phi_i(R_i + \gamma\Delta_i)\\
        & \text{subject to} && \Delta_i \ge 0,
    \end{aligned}
\]
and the second case has $\Delta_i = 0$ and $\Lambda_i \ge 0$, which means
that we only have to solve
\[
    \begin{aligned}
        & \text{maximize} && \eta_i \Lambda_i + \lambda_i \phi_i(R_i - \Lambda_i)\\
        & \text{subject to} && \Lambda_i \ge 0.
    \end{aligned}
\]
It would then suffice to take whichever of the two cases has the highest
optimal objective value---though, unless $\gamma = 1$, at most one problem will
have a positive solution and we deal with the $\gamma = 1$ case below. These
problems can be solved by ternary search (if we only have access to $\phi_i$
via function evaluations), bisection (if we also have access to the derivative,
$\phi_i'$), or Newton's method (if we have access to the second derivative,
$\phi_i''$).

These problems also often have closed form solutions. For example, the
optimality conditions for the first of the two cases is: if $\Delta^\star_i =
0$ is optimal, then
\[
    \lambda_i \gamma\phi_i'(R_i) \le \eta_i,
\]
or, otherwise, $\Delta_i^\star > 0$ satisfies
\[
    \lambda_i \gamma\phi_i'(R_i + \gamma\Delta_i^\star) = \eta_i.
\]
The former condition is a simple check and the latter condition is a simple
root-finding problem that, in many cases, has a closed-form solution. A very
similar also holds for the second case.

Finally, if $\gamma = 1$, the subproblems in the dual
function~\eqref{eq:app-dual-sep} simplify even further to the unconstrained
single variable convex optimization problem
\[
    \sup_{t} \left(\eta_i t + \lambda_i \phi_i(R_i - t)\right),
\]
which is easily solved via any number of methods. We can recover a solution to
the original subproblem by setting $\Lambda_i^\star =  t^\star +
\Delta_i^\star$ for any solution $t^\star$, where $\Delta_i^\star$ is any
value.

\paragraph{CFMMs as conic constraints.}
Define the power cone as 
\[
    K_\mathrm{pow}(w) = \{ (x,y,z) \in \reals^3 \mid 
  x^{w} y^{1-w} \ge |z|, ~x \ge 0, ~y \ge 0 \}.
\]
We model the two-asset market constraints as
\begin{equation}\label{eq:mkt}
    (R + \gamma \Delta - \Lambda, \,\phi(R)) \in K_\mathrm{pow}(w),
    \qquad \text{and} \qquad 
    \Delta, \Lambda \ge 0,
\end{equation}
where $w = 0.5$ for Uniswap and $w = 0.8$ for Balancer.
Define the geometric mean cone as 
\[
    K_\mathrm{geomean} = \left\{ (t, x) \in
      \reals \times \reals^{n} \mid x \ge 0,~ \left(x_1x_2\cdots x_n\right)^{1/n} \ge t \right\}
\]
We model the multi-asset market constraint as
\begin{equation}
    (-3\phi(R), \, R + \gamma \Delta - \Lambda) \in K_\mathrm{geomean},
    \qquad \text{and} \qquad 
    \Delta, \Lambda \ge 0.
\end{equation}

\paragraph{Objectives as conic constraints.}
Define the rotated second order cone as
\[
    K_\mathrm{rot2} = \{ (t,u,x) \in \reals_+ \times \reals_+ \times \reals^{n} \mid
     2tu \ge \| x \|_2^2\}.
\]
The net flow utility function is
\[
    U(y) = c^Ty - (1/2)\sum_{i=1}^n (y_i)_-^2,
\]
where $x_- = \max(-x, 0)$ denotes the negative part of $x$. In conic form,
maximizing $U$ is equivalent to maximizing
\[
    c^Ty - (1/2)\sum_{i=1}^n (p_1)_i
\]
subject to the constraints
\[
    p_2 \ge 0, \qquad p_2 \ge -y, \qquad (p_1)_i, (p_2)_i) \in K_\mathrm{rot2} \quad \text{for} ~ i = 1, \dots, n,
\]
where we introduced new variables $p_1, p_2 \in \reals^n$~\cite[\S3.2.2]{mosekcookbook}. 
The $V_i$'s can be modeled similarly using the rotated second order cone.

\paragraph{Conic form problem.}
The CFMM arbitrage example can then be written in conic form as
\[
\begin{aligned}
    & \text{maximize} && c^Ty - (1/2)\sum_{i=1}^n (p_1)_i  - (1/2)\sum_{i=1}^m (t_1)_i \\
    & \text{subject to} && (0.5,\, (p_1)_i,\, (p_2)_i) \in K_\mathrm{rot2}, \quad i = 1, \dots, n \\
    &&& p_1 \ge 0 \\
    &&&  p_2 \ge 0, \quad p_2 \ge -y \\
    &&& (0.5,\, (t_1)_i,\, (t_2)_i) \in K_\mathrm{rot2}, \quad i = 1, \dots, n \\
    &&& t_1 \ge 0 \\
    &&& t_2 \ge 0, \quad (t_2)_i \ge -(\Lambda_i - \Delta_i) \\
    &&& (R + \gamma \Delta - \Lambda,\,  \phi(R)) \in K_\mathrm{pow}(w_i), \quad i = 1, \dots, m_1 \\
    &&& (-3\phi(R),\, R + \gamma \Delta - \Lambda) \in K_\mathrm{geomean}, \quad i = m_1 + 1, \dots, m \\\
    &&& \Delta_i,\; \Lambda_i \ge 0, \quad i = 1, \dots, m,
\end{aligned}
\]
with variables $y \in \reals^n$, $p_1 \in \reals^n$, $p_2 \in \reals^n$, $t_1 \in \reals^m$, 
$(t_2)_i \in \reals^{n_i}$, $\Delta \in \reals^{n_i}$, and $\Lambda \in \reals^{n_i}$ for $i = 1, \dots, m$.

\section{Automated conservation laws}\label{app:conservation} 
We define the set of conservation laws for an instance of the convex flow
problem as
\[
    C = \left\{
        c \in \reals^n 
        ~\middle|~
        c^T \left(\sum_{i=1}^m A_i x_i\right) \ge 0 
        ~~\text{for all}~~ 
        x_i \in T_i
    \right\}.
\]
This construction is a generalization of the conservation law discussed
in~\S\ref{sec:app-max-flow}. Our goal is to find vectors in this set. Define the
dual cone for a set $S$ as
\[
    K(S) = \{ y \mid y^Tx \ge 0 ~~\text{for all}~~ x \in S \}.
\]
A sufficient condition for a vector $c$ to be in $C$ is that 
\[
    c \in \bigcap_{i=1}^m K(A_i T_i).
\]
(Of course, in general, this condition is not necessary.) We can then find
conservation laws by solving the convex optimization problem
\[
\begin{aligned}
    &\text{find}        && c \\
    &\text{subject to}  && c^Tz = 1 \\
    &                   && c \in \bigcap_{i=1}^m K(A_i T_i),
\end{aligned}
\]
where $z$ is some nonzero vector, for example, sampled from the standard normal
distribution.

\end{document}